\documentclass[11pt]{amsart}
\usepackage[all]{xy}

\usepackage{url}

\numberwithin{equation}{section}
\numberwithin{figure}{section}

\newtheorem{Thm}{Theorem}
\newtheorem{Lem}{Lemma}[subsection]
\newtheorem{Cor}[Lem]{Corollary}
\newtheorem{Prop}[Lem]{Proposition}

\newcommand{\Aa}{\mathsf{A}}
\newcommand{\Dd}{\mathsf{D}}
\newcommand{\Ee}{\mathsf{E}}

\newcommand{\Z}{\mathbb{Z}}     % integers
\newcommand{\N}{\mathbb{N}}     % naturals (with 0)
\newcommand{\C}{\mathbb{C}}     % complex numbers

\newcommand{\PP}{\mathbb{P}}

\newcommand{\bQ}{\overline{Q}}

\newcommand{\bil}[1]{\langle #1\rangle}
\newcommand{\abs}[1]{\left| #1\right|}
\newcommand{\str}[1]{\langle #1\rangle} % \del-stratification
\newcommand{\mt}[1]{\operatorname{mt}(#1)} % middle term of exact sequence

\newcommand{\id}{1\kern -.35em 1}

\newcommand{\md}{\operatorname{mod}}
\newcommand{\smd}{\underline{\md}}

\newcommand{\Hom}{\operatorname{Hom}}
\newcommand{\sthom}{\operatorname{\underline{Hom}}}
\newcommand{\Ext}{\operatorname{Ext}}

\newcommand{\dimv}{\operatorname{\underline{\dim}}}
\newcommand{\Gl}{\operatorname{GL}}
\newcommand{\Tr}{\operatorname{Tr}}

\newcommand{\Ima}{\operatorname{Im}}
\newcommand{\Ker}{\operatorname{Ker}}

\newcommand{\ka}{\mathrm{k}}

\newcommand{\La}{\Lambda}
\newcommand{\Ga}{\Gamma}
\newcommand{\Sig}{\Sigma}
\newcommand{\Ome}{\Omega}

\newcommand{\alp}{\alpha}
\newcommand{\bet}{\beta}
\newcommand{\eps}{\epsilon}
\newcommand{\del}{\delta}
\newcommand{\lam}{\lambda}
\newcommand{\vph}{\varphi}

\newcommand{\De}{{\mathcal D}}
\newcommand{\M}{{\mathcal M}}
\newcommand{\SB}{{\mathcal S}}
\newcommand{\cT}{{\mathcal T}}

\newcommand{\cc}{\mathbf{c}}
\newcommand{\dd}{\mathbf{d}}
\newcommand{\ee}{\mathbf{e}}
\newcommand{\ii}{\mathbf{i}}
\newcommand{\jj}{\mathbf{j}}

\def\V{{\mathbf V}}
\def\U{{\mathbf U}}
\def\W{{\mathbf W}}
\def\T{{\mathbf T}}

\newcommand{\f}{\mathfrak{f}}
\newcommand{\g}{\mathfrak{g}}

\newcommand{\n}{\mathfrak{n}}

\def\ra{\rightarrow}

\newcommand{\xra}[1]{\xrightarrow{#1}}

\newcommand{\ol}[1]{\overline{#1}}
\newcommand{\mc}{\text{-}}

\newcommand{\df}{\colon}

\begin{document}

%\date{January 2007}

\bigskip
\title[Semicanonical bases and preprojective algebras II]%
{Semicanonical bases and preprojective algebras II:\\ A multiplication
formula}

\author{Christof Gei{\ss}}
\address{
Christof Gei{\ss}\newline
Instituto de Matem\'aticas\newline 
Universidad Nacional Aut\'onoma de M\'exico\newline 
Ciudad Universitaria\newline
04510 M\'exico D.F.\newline 
M\'exico}
\email{christof@math.unam.mx}
\urladdr{\url{www.matem.unam.mx/~christof/}}

\author{Bernard Leclerc}
\address{
Bernard Leclerc\newline
Laboratoire LMNO\newline 
Universit\'e de Caen\newline 
F-14032 Caen Cedex\newline 
France}
\email{leclerc@math.unicaen.fr}
\urladdr{\url{www.math.unicaen.fr/~leclerc/}}

\author{Jan Schr\"oer}
\address{
Jan Schr\"oer\newline
Mathematisches Institut\newline 
Universit\"at Bonn\newline 
Beringstr. 1\newline
53115 Bonn\newline 
Germany}
\email{schroer@math.uni-bonn.de}
\urladdr{\url{http://www.math.uni-bonn.de/people/schroer/}}

\subjclass[2000]{14M99, 16G20, 17B35, 17B67, 20G05}
\keywords{Universal enveloping algebra, semicanonical basis, preprojective
algebra, flag variety, composition series}

%%%%%%%%%%%%%%%%%%%%%%%%%%%%%%%%%%%%%%%%%%%%%%%%%%%%%%%%%%%%%%%%%%%%%%%

\begin{abstract}
Let $\n$ be a maximal nilpotent subalgebra of a complex symmetric Kac-Moody Lie
algebra. 
Lusztig has introduced a basis of $U(\n)$ called the
semicanonical basis, whose elements can be seen as certain constructible
functions on varieties of nilpotent modules over a preprojective algebra of the
same type as $\n$.
We prove a formula for the product of two elements of the
dual of this semicanonical basis, and more generally  
for the product of two evaluation forms
associated to arbitrary modules over the preprojective algebra.
This formula plays an important role in our work on the relationship 
between semicanonical bases, representation theory 
of preprojective algebras, and Fomin and Zelevinsky's theory of cluster
algebras. It was inspired by recent results of Caldero and Keller. 
\end{abstract}

\maketitle

%%%%%%%%%%%%%%%%%%%%%%%%%%%%%%%%%%%%%%%%%%%%%%%%%%%%%%%%%%%%%%%%%%%%%%%

\section{Introduction and main result}

%%%%%%%%%%%%%%%%%%%%%%%%%%%%%%%%%%%%%%%%%%%%%%%%%%%%%%%%%%%%%%%%%%%%%%%

%%%%%%%%%%%%%%%%%%%%%%%%%%%%%%%%
\subsection{Semicanonical bases}
%%%%%%%%%%%%%%%%%%%%%%%%%%%%%%%%
Let $\La$ be the preprojective algebra associated to a 
connected quiver $Q$ without loops. This is the  associative algebra 
\[
\La = \C\bQ /\bil{\sum_{a\in Q_1}(\bar{a}a-a\bar{a})},
\]
where $\bQ$ denotes the double of $Q$ and $Q_1$ is the set of 
arrows of $Q$ (see {\em e.g.}~\cite{Ri1,GLS}). 
Recall that $\bQ$ is obtained from $Q$ by inserting for each arrow 
$a\df i\ra j$ in $Q$ a new arrow $\bar{a}\df j\ra i$.
The algebra
$\La$ is independent of the orientation of $Q$, and it is finite-dimensional
(and selfinjective) if and only if the underlying graph $\Delta$ of $Q$ is
a simply laced Dynkin diagram, {\em i.e.} if
$\Delta\in\{ \Aa_n (n \ge 1), \, \Dd_n (n \ge 4), \, \Ee_n (n = 6,7,8) \}$.
In the sequel, all $\La$-modules are assumed to be finite-dimensional 
nilpotent left modules.
We denote by $I$ the set of vertices of $Q$, and by $\La_\dd$ the
affine variety of nilpotent  
$\La$-modules with dimension
vector $\dd=(d_i)_{i\in I}$. The reductive group
$\Gl_{\dd}=\prod_{i\in I} \Gl_{d_i}(\C)$
acts on $\La_\dd$ by base change.

Let $\n$ be a maximal nilpotent subalgebra of a
complex Kac-Moody Lie algebra $\g$ of type $\Delta$.
Lusztig \cite{Lu,Lu2} proved that the enveloping algebra $U(\n)$ is 
isomorphic to 
\[
\M = \bigoplus_{\dd\in\N^I} \M(\dd), 
\]
where the $\M(\dd)$ are certain vector spaces of 
$\Gl_{\dd}$-invariant constructible functions 
on the affine varieties $\La_\dd$ 
(see~\ref{ssec:repeat} below for the definition of $\M(\dd)$).
This yields a new basis $\SB$ of $U(\n)$
indexed by the irreducible components of the varieties $\La_\dd$, called
the {\em semicanonical basis} \cite{Lu2}.

Let $\M^*$ be the graded dual of $\M$. 
A multiplication on $\M^*$ is defined via the natural comultiplication
of the Hopf algebra $U(\n) \cong\M$.
In \cite{GLS} we considered the basis $\SB^*$ of $\M^*$ dual to the
semicanonical basis of $\M$, and began to study its multiplicative properties.

%%%%%%%%%%%%%%%%%%%%%%%%%%%%%%%%%%%%%%%%%%%%%%%%%%%%%%
\subsection{$\del$-Stratification} \label{deltafinite} 
%%%%%%%%%%%%%%%%%%%%%%%%%%%%%%%%%%%%%%%%%%%%%%%%%%%%%%
For a $\La$-module $x \in \La_\dd$ define the evaluation form 
$
\delta_x\df \M \to \C 
$
which maps a constructible function $f \in \M(\dd)$ to $f(x)$.
Define
\[
\str{x} := \{y\in\La_\dd\mid \del_y=\del_x\}.
\]
This is a constructible subset of $\La_\dd$ since $\M(\dd)$ is a 
finite-dimensional space of constructible functions on $\La_\dd$. 
For the same reason
we can choose a {\em finite} set $R(\dd)\subset\La_\dd$ such that
\[
\La_\dd = \bigsqcup_{x\in R(\dd)} \str{x}.
\]
Each irreducible component $Z$ of $\La_\dd$ has a unique stratum
$\str{x}\cap Z$ containing a dense open
subset of $Z$, and the points of this stratum are called the {\em generic} points of $Z$.
We can then reformulate the definition of $\SB^*$ as follows:
the element $\rho_Z$ of $\SB^*$ labeled by $Z$ is equal to
$\del_x$ for a generic point $x$ of $Z$. 

The aim of this paper is to prove a general formula for the product of two
arbitrary evaluation forms $\del_x$.

%%%%%%%%%%%%%%%%%%%%%%%%%%%%%%%%%%%%%%%%%
\subsection{Extensions}\label{extensions}
%%%%%%%%%%%%%%%%%%%%%%%%%%%%%%%%%%%%%%%%% 
Let $x'\in\La_{\dd'}$, $x''\in\La_{\dd''}$, and $\ee:=\dd'+\dd''$. 
For a constructible 
$\Gl_{\ee}$-invariant subset $S\subseteq \La_{\ee}$ we consider 
\[
\Ext^1_\La(x',x'')_S=\{[ 0\ra x''\ra y\ra x'\ra 0]\in
\Ext^1_\La(x',x'')\setminus\{0\}\mid y\in S\}.
\]
This is a constructible $\C^*$-invariant subset of the quasi-affine space
$\Ext^1_\La(x',x'')\setminus\{0\}$, see \ref{E-del}. Thus we may consider the 
constructible subset $\PP\Ext^1_\La(x',x'')_S := \Ext^1_\La(x',x'')_S/\C^*$
of the projective space $\PP\Ext^1_\La(x',x'')$. 
For a constructible subset $U$ of a complex variety let $\chi(U)$ be its
(topological) Euler characteristic with respect to cohomology with compact 
support.
For the empty set $\emptyset$ we set $\chi(\emptyset) = 0$.
Note that we have, by additivity of $\chi$,
\[
\sum_{x\in R(\ee)} \chi(\PP\Ext^1_\La(x',x'')_{\str{x}}) =
\chi(\PP\Ext^1_\La(x',x''))=
\dim \Ext^1_\La(x',x''),
\]
see also~\ref{pushouts}.

%%%%%%%%%%%%%%%%%%%%%%%%
\subsection{Main result}
%%%%%%%%%%%%%%%%%%%%%%%%

\begin{Thm}[Multiplication Formula] \label{thm1}
With the notation of \ref{extensions}, we have
\[
\chi(\PP\Ext^1_\La(x',x''))\, \del_{x'\oplus x''} =
\sum_{x\in R(\ee)}
\left(\chi(\PP\Ext^1_\La(x',x'')_{\str{x}})+
\chi(\PP\Ext^1_\La(x'',x')_{\str{x}})\right)\,\del_x.
\]
\end{Thm}

This theorem is very much inspired by recent work 
by Caldero and Keller on cluster algebras and cluster categories of finite 
type \cite{CK}.

Note that our multiplication formula is meaningless if $\Ext^1_\La(x',x'')=0$
(or equivalently if $\Ext^1_\La(x'',x')=0$). 
We proved in \cite{GLS} that
in all cases 
$$
\del_{x'\oplus x''}=\del_{x'}\cdot\del_{x''}.
$$ 
For this reason we regard Theorem~\ref{thm1} as a multiplication formula. 
Using the  notation from~\cite{GLS} the formula of the theorem can be 
rewritten as
\[
\del_{x'}\cdot\del_{x''}=\frac{1}{\dim \Ext^1_\La(x',x'')}
\left(\int_{\C^*[\eta]\in\PP\Ext^1_\La(x',x'')}\del_{\mt{\eta}} +
      \int_{\C^*[\eta]\in\PP\Ext^1_\La(x'',x')}\del_{\mt{\eta}}\right).
\]
Here, the ``integrals'' are taken with respect to the ``measure'' given
by the Euler characteristic, $\mt{\eta}$ stands for the
isomorphism class of the middle term of the extension $\eta$.
If $[\eta]\neq 0$ we write $\C^*[\eta]$ for the corresponding element in the 
projective space $\PP\Ext^1_\La(x',x'')$.

The following important special case is easier to state and also
easier to prove. 

\begin{Thm} \label{mainresult}
Let $x'$ and $x''$ be $\La$-modules such that $\dim \Ext_\La^1(x',x'') = 1$,
and let 
\[
0 \to x'' \to x \to x' \to 0
\text{\;\;\; and \;\;\;}
0 \to x' \to y \to x'' \to 0
\]
be non-split short exact sequences.
Then
$\delta_{x'} \cdot \delta_{x''} = \delta_{x} + \delta_{y}$.
\end{Thm}

A $\La$-module $x$ is called {\em rigid} if $\Ext_\La^1(x,x) = 0$.
If $x$ is rigid, it is generic and $\delta_x$ is a dual semicanonical 
basis vector.
In the above theorem, if $x'$ and $x''$ are rigid modules with 
$\dim \Ext_\La^1(x',x'') = 1$ one can show that $x$ and $y$ have to be rigid as
well, see~\cite{GLS3}.
Thus we obtain a multiplication formula for certain dual semicanonical
basis vectors.
This is analyzed and interpreted in more detail in~\cite{GLS3}.

%%%%%%%%%%%%%%%%%%%%%%%%%%%%%%%%%%%%%%%%%%%%%%
\subsection{Calabi-Yau property}\label{ssecCY}
%%%%%%%%%%%%%%%%%%%%%%%%%%%%%%%%%%%%%%%%%%%%%%
For  the proof of our result it is crucial that in the category
of finite-dimensional (nilpotent) $\La$-modules we have a functorial pairing
$$
\Ext^1_\La(x,y)\times \Ext^1_\La(y,x)\ra\C.
$$ 
In other words, for
each pair $(x,y)$ of finite-dimensional (nilpotent) 
$\La$-modules there is an isomorphism
\[
\phi_{x,y}\df \Ext^1_\La(x,y)\ra D\Ext^1_\La(y,x)
\]
such that for any $\lam\in\Hom_\La(y,y')$, $[\eta]\in\Ext^1_\La(x,y)$,
$\rho\in\Hom_\La(x',x)$ and $[\eps]\in\Ext^1_\La(y',x')$ one has
\begin{equation} \label{eqn:CY-Func}
\phi_{x',y'}(\lam\circ [\eta]\circ\rho)([\eps])=
\phi_{x,y}([\eta])(\rho\circ[\eps]\circ\lam).
\end{equation}
This is well-known to specialists. The usual argument is that this
is the ``shadow'' of certain 2-Calabi-Yau properties for preprojective
algebras, see \ref{rec:2cy} for more details.

However, in the non-Dynkin case this argument relies on the fact that 
then the preprojective algebra is a 2-Calabi-Yau algebra. 
Unfortunately, we could not localize a written proof for this last
fact. So we include in Section~\ref{sec:extsym} a direct proof for
the required functorial isomorphism which does not depend on the type
of $Q$ and neither requires nilpotency for the finite-dimensional modules.

%%%%%%%%%%%%%%%%%%%%%%%%%%%%%%%%%%%%%%%%%%%%%%%%%%%%%%%%%%%%%%%%%%%%%%%

\section{Flags and composition series induced by short exact sequences}

%%%%%%%%%%%%%%%%%%%%%%%%%%%%%%%%%%%%%%%%%%%%%%%%%%%%%%%%%%%%%%%%%%%%%%%

%%%%%%%%%%%%%%%%%%%%%%%%%%%%%%%%%%%%%%%%%%%%%%%%%%%%%%%%%
\subsection{Definition of $\M(\dd)$.} \label{ssec:repeat}
%%%%%%%%%%%%%%%%%%%%%%%%%%%%%%%%%%%%%%%%%%%%%%%%%%%%%%%%%
We repeat some notation from \cite{GLS}.
Recall that $I$ is the set of vertices of $Q$.
We identify the elements of $I$ with the natural basis of $\Z^I$. Thus, for
$\ii = (i_1,\ldots ,i_m)$  a sequence of elements of $I$, we may define
$\abs{\ii} := \sum_{k=1}^m i_k\in\Z^I$.

Let $\V$ be an $I$-graded vector space with graded dimension
$|\V|\in\Z^I$.
We denote by $\La_\V$ the affine variety of $\La$-module structures
on $\V$. Clearly, if $\V'$ is another graded space with $|\V'|=|\V|$
the varieties $\La_\V$ and $\La_{\V'}$ are isomorphic.
Therefore, it is sometimes convenient to write
$\La_{|\V|}$ instead of $\La_\V$ and to think of it as the variety
of $\La$-modules with dimension vector $\dd=|\V|$.

If $\cc = (c_1,\ldots ,c_m)\in \{0,1\}^m$ and $\V$ is
an $I$-graded vector space such that $\sum_k c_ki_k = |\V|$, we 
define a {\em flag} in $\V$ of type $(\ii,\cc)$ as a sequence
\[
\f = \left(\V=\V^0 \supseteq \V^1 \supseteq \cdots \supseteq \V^m = 0\right)
\]
of $I$-graded subspaces of $\V$ such that 
$
|\V^{k-1}/\V^k| = c_k i_k
$
for $k=1,\ldots, m$.
We denote by $\Phi_{\ii,\cc}$ the variety of flags of type
$(\ii,\cc)$.
When $(c_1,\ldots ,c_m) = (1,\ldots ,1)$, we simply write $\Phi_{\ii}$.

Let $x \in \La_\V$. 
A flag $\f$ is $x$-{\it stable} 
if $x(\V^k)\subseteq \V^k$ for all $k$.
We denote by $\Phi_{\ii,\cc,x}$ (resp. $\Phi_{\ii,x}$) 
the variety of $x$-stable flags of type
$(\ii,\cc)$ (resp. of type $\ii$).
Note that an $x$-stable flag is the same as a composition series
of $x$ regarded as a $\La$-module.
Let
$
d_{\ii,\cc}\df \La_\V \to \C 
$
be the  function defined by 
$$
d_{\ii,\cc}(x) = \chi(\Phi_{\ii,\cc,x}).
$$
It follows from \cite[\S 5.8]{GLS} that $d_{\ii,\cc}$ is a constructible
function, see also~\cite[\S 2.1]{Lu}.
If $c_k=1$ for all $k$, we simply write $d_\ii$ instead of $d_{\ii,\cc}$.
In general, $d_{\ii,\cc}=d_\jj$ where $\jj$ is the subword of $\ii$
consisting of the letters $i_k$ for which $c_k=1$.
We then define $\M(\V)$ as the vector space over $\C$ spanned by
the functions $d_\jj$, where $\jj$ runs over
all words in $I$ with $\abs{\jj}=|\V|$.
Again this space only depends on $\dd=|\V|$, up to isomorphism,
and we also denote it by $\M(\dd)$.
This is equivalent to Lusztig's definition in \cite{Lu, Lu2}. 

%%%%%%%%%%%%%%%%%%%%%%%%%%%%%%%%%%%%%%%%%%%%%%%
\subsection{Reformulation}\label{reformulation}
%%%%%%%%%%%%%%%%%%%%%%%%%%%%%%%%%%%%%%%%%%%%%%%
In order to prove Theorem \ref{thm1} and
Theorem \ref{mainresult} it is sufficient to show that both sides of
the respective equalities coincide after evaluation at an arbitrary
$d_{\ii}$, see~\ref{ssec:repeat}.
Since 
$$
\delta_x(d_\ii)=\chi(\Phi_{\ii,x})
$$ 
to prove Theorem \ref{thm1} we have to show that 
\begin{equation}\label{eqn:ref1}
\chi(\PP\Ext^1_\La(x',x''))\cdot \chi(\Phi_{\ii,x'\oplus x''})=
\sum_{x\in R(\ee)} 
\left(\chi(\PP\Ext^1_\La(x',x'')_{\str{x}})+
\chi(\PP\Ext^1_\La(x'',x')_{\str{x}})\right) \cdot \chi(\Phi_{\ii,x})
\end{equation}
for all sequences $\ii$ with $\abs{\ii}=\ee=\dimv(x')+\dimv(x'')$.
Similarly, for Theorem \ref{mainresult} we have to show that
$$
\chi(\Phi_{\ii,x'\oplus x''})=\chi(\Phi_{\ii,x})+\chi(\Phi_{\ii,y})
$$
for all sequences $\ii$ with $\abs{\ii}=\dimv(x')+\dimv(x'')$.  

%%%%%%%%%%%%%
\subsection{}\label{ssec:refine}
%%%%%%%%%%%%%
Let $\V',\V'',\V$ be $I$-graded vector spaces such that
$|\V'| + |\V''| = |\V|$, and
let $\ii=(i_1,\ldots,i_m) \in I^m$ with $\sum_k i_k = |\V|$.
Let $x'\in\La_{\V'}$, $x''\in\La_{\V''}$ 
and $x \in \La_\V$.
For a short exact sequence 
\[
\eps: \; 0 \to x'' \xrightarrow{\iota} x \to x' \to 0
\]
define a map
\[
\alpha_{\ii,\eps}\df  \Phi_{\ii,x} \to 
\bigsqcup_{(\cc',\cc'')} \Phi_{\ii,\cc',x'} \times 
\Phi_{\ii,\cc'',x''}
\]
which maps a flag
$
\f_x = (\V^l)_{0\le l\le m}\in\Phi_{\ii,x}
$
of submodules of $x$ 
to 
$(\f_{x'},\f_{x''})$ where 
\[
\f_{x'}=(\V^l/(\V^l \cap \iota(x'')))_{0\le l\le m}
\text{\;\;\; and \;\;\;}
\f_{x''}=(\V^l \cap \iota(x''))_{0\le l\le m}\,.
\]
Here we regard $\f_{x''}$ as a flag in $\V''$ by identifying $\V''$ with
$\iota(x'')$,
and $\f_{x'}$ as
a flag in $\V'$ by identifying $\V'$ with $\V/\iota(x'')$.
Clearly, we have 
$(\f_{x'},\f_{x''}) \in \Phi_{\ii,\cc',x'} \times \Phi_{\ii,\cc'',x''}$ 
for some sequences $\cc'$, $\cc''$ in $\{0,1\}^m$ satisfying 
\begin{equation}\label{condition}
\sum_{k=1}^m c'_ki_k=|\V'|,
\qquad \sum_{k=1}^m c''_ki_k=|\V''|,
\qquad c'_k + c''_k = 1\ (1 \le k \le m).
\end{equation} 
Let $W_{\V',\V''}$ denote the set of pairs $(\cc',\cc'')$ satisfying
(\ref{condition}).
Define
\begin{align*}
L^1_{\ii,\cc',\cc'',\eps} &= 
(\Phi_{\ii,\cc',x'} \times \Phi_{\ii,\cc'',x''}) \cap 
\Ima(\alpha_{\ii,\eps}),\\
L^2_{\ii,\cc',\cc'',\eps} &= 
(\Phi_{\ii,\cc',x'} \times \Phi_{\ii,\cc'',x''}) \setminus 
L^1_{\ii,\cc',\cc'',\eps},\\
\Phi_{\ii,x}(\cc',\cc'',\eps) &=
\alpha_{\ii,\eps}^{-1}(\Phi_{\ii,\cc',x'} \times \Phi_{\ii,\cc'',x''}).
\end{align*}
Thus we obtain maps 
\[
\alpha_{\ii,\eps}(\cc',\cc'')\df \Phi_{\ii,x}(\cc',\cc'',\eps)
\to \Phi_{\ii,\cc',x'} \times \Phi_{\ii,\cc'',x''}.
\]
Set
$W_{\V',\V'',\eps} = \{ (\cc',\cc'') \in W_{\V',\V''} \mid 
\Phi_{\ii,x}(\cc',\cc'',\eps) 
\not= \emptyset \}$.
We get a finite partition
\[
\Phi_{\ii,x}
= \bigsqcup_{(\cc',\cc'') \in W_{\V',\V'',\eps}}
\Phi_{\ii,x}(\cc',\cc'',\eps).
\]

%%%%%%%%%%%%%
\subsection{}
%%%%%%%%%%%%%
For the rest of this section let
\[
\f_{x'} = 
(x' = x_0' \supseteq x_1' \supseteq  \cdots \supseteq 
x_m' = 0)
\text{\;\;\; and \;\;\;}
\f_{x''} = 
(x'' = x_0'' \supseteq x_1'' \supseteq \cdots \supseteq 
x_m'' = 0)
\]
be flags 
with $(\f_{x'},\f_{x''}) \in \Phi_{\ii,\cc',x'} \times \Phi_{\ii,\cc'',x''}$
for some $(\cc',\cc'') \in W_{\V',\V''}$ where 
$x' \in \La_{\V'}$ and $x'' \in \La_{\V''}$.

For $1 \le j \le m$ let $\iota_{x',j}$ and $\iota_{x'',j}$ be
the inclusion maps $x_j' \to x_{j-1}'$ and
$x_j'' \to x_{j-1}''$, respectively.
Note that for all $j$ either $\iota_{x',j}$ or $\iota_{x'',j}$ is an identity
map.
In particular, either $x'_{m-1} = 0$ or $x''_{m-1} = 0$.

All results in this section and also the definition of the maps
$\beta'$ and $\beta$ (see below) are inspired by \cite{CK}.
The following lemma follows directly from the definitions and the 
considerations in Section \ref{pushouts}.

\begin{Lem}\label{imagealpha}
Let 
$
\eps: \; 0 \to x'' \to x \to x' \to 0
$
be a short exact sequence of $\La$-modules.
Then the following are equivalent:
\begin{itemize}

\item[(i)]
$(\f_{x'},\f_{x''}) \in L^1_{\ii,\cc',\cc'',\eps}$;

\item[(ii)]
There exists a commutative diagram with exact rows of the form
\[
\xymatrix{
\eps: & 0 \ar[r] & x'' \ar[r] & x \ar[r] & x' \ar[r] & 0\\ 
\eps_0: & 0 \ar[r]  & x_0'' \ar[r]\ar@{=}[u] & x_0 \ar[r]\ar@{=}[u] & 
x_0' \ar[r]\ar@{=}[u] & 0\\
\eps_1: & 0 \ar[r]  & x''_1 \ar[r]\ar[u]^{\iota_{x'',1}} & 
x_1 \ar[r] \ar[u] & 
x'_1 \ar[r]\ar[u]_{\iota_{x',1}} & 0\\
& & \vdots \ar[u]^{\iota_{x'',2}}& \vdots \ar[u]& \vdots 
\ar[u]_{\iota_{x',2}}\\
\eps_{m-2}: & 0 \ar[r]  & x''_{m-2} \ar[r]\ar[u]^{\iota_{x'',m-2}} & 
x_{m-2} \ar[r] \ar[u] & 
x'_{m-2} \ar[r]\ar[u]_{\iota_{x',m-2}} & 0\\
\eps_{m-1}: & 0 \ar[r] & x''_{m-1} \ar[r]\ar[u]^{\iota_{x'',m-1}} & 
x_{m-1} \ar[r]\ar[u] & 
x'_{m-1} \ar[r]\ar[u]_{\iota_{x',m-1}} & 0
}
\]
\item[(iii)]
There exist elements $[\eps_i]\in\Ext^1_\La(x'_i,x''_i)$ for 
$i=0,1,\ldots, m-1$ such that $[\eps_0]=[\eps]$ and
\[
[\eps_{j-1}]\circ \iota_{x',j}= \iota_{x'',j}\circ [\eps_j]
\text{ for } j=1,2,\ldots, m-1.
\]
\end{itemize}
\end{Lem}

We define a map 
\begin{multline} \label{eqn:betp}
\beta'_{\ii,\cc',\cc'',\f_{x'},\f_{x''}}\df 
\bigoplus_{j=0}^{m-2} \Ext^1_\La(x_j',x_j'')=:W' \to 
V':=\bigoplus_{j=0}^{m-2} \Ext^1_\La(x_{j+1}',x_j'') \\
([\eps_0],\ldots,[\eps_{m-2}]) \mapsto 
\sum_{j=0}^{m-3}([\eps_j]\circ\iota_{x',j+1}-\iota_{x'',j+1}\circ[\eps_{j+1}])+
[\eps_{m-2}]\circ\iota_{x',m-1}
\end{multline}

Observe that
$[\eps_j]\circ\iota_{x',j+1}-\iota_{x'',j+1}\circ[\eps_{j+1}]$ 
is an element of $\Ext_\La^1(x_{j+1}',x_j'')$.
By
$$
p_0\df \bigoplus_{j=0}^{m-2} \Ext_\La^1(x_j',x_j'') \to 
\Ext_\La^1(x_0',x_0'')
$$
we denote the canonical projection map.

\begin{Lem}\label{betadash}
For a short exact sequence
$
\eps\df 0 \to x'' \to x \to x' \to 0
$,
the following are equivalent:
\begin{itemize}

\item[(i)]
$(\f_{x'},\f_{x''}) \in L^1_{\ii,\cc',\cc'',\eps}$;

\item[(ii)]
$[\eps] \in p_0(\Ker(\beta'_{\ii,\cc',\cc'',\f_{x'},\f_{x''}}))$.

\end{itemize}
\end{Lem}

\begin{proof}
We have $[\eps] \in p_0(\Ker(\beta'_{\ii,\cc',\cc'',\f_{x'},\f_{x''}}))$
if and only if 
$
\beta'_{\ii,\cc',\cc'',\f_{x'},\f_{x''}}([\eps_0],\ldots,[\eps_{m-2}]) 
= 0
$
for some $([\eps_0],\ldots,[\eps_{m-2}])$ 
with $[\eps_0] = [\eps]$.
This is the case if and only if 
$\iota_{x'',j+1}\circ[\eps_{j+1}]=[\eps_j]\circ\iota_{x',j+1}$
for all $0 \le j \le m-3$, $[\eps_{m-2}]\circ\iota_{x',m-1} = 0$ and 
$[\eps_0] = [\eps]$.
Now the result follows from Lemma \ref{imagealpha}.
\end{proof}

Next, we define a map 
\begin{multline} \label{eqn:bet}
\beta_{\ii,\cc'',\cc',\f_{x''},\f_{x'}}\df 
\bigoplus_{j=0}^{m-2} \Ext^1_\La(x_j'',x_{j+1}')=: V \to 
W:=\bigoplus_{j=0}^{m-2} \Ext^1_\La(x_j'',x_j')\\
([\eta_0],\ldots,[\eta_{m-2}]) \mapsto 
\iota_{x',1}\circ[\eta_0] + 
\sum_{j=1}^{m-2} (\iota_{x',j+1}\circ [\eta_j] - [\eta_{j-1}]\circ\iota_{x'',j})
\end{multline}

Note that 
$\iota_{x',j+1}\circ[\eta_j] - [\eta_{j-1}]\circ\iota_{x'',j}$ 
is an element of $\Ext_\La^1(x_j'',x_j')$.

\begin{Lem}\label{beta}
For a short exact sequence
$
\eta\df 0 \to x' \to y \to x'' \to 0
$
the following are equivalent:
\begin{itemize}
\item[(i)]
$(\f_{x''},\f_{x'}) \in L^1_{\ii,\cc'',\cc',\eta}$;
\item[(ii)]
$[\eta] \in \Ext_\La^1(x'',x') \cap 
\Ima(\beta_{\ii,\cc'',\cc',\f_{x''},\f_{x'}})$.
\end{itemize}
\end{Lem}

\begin{proof}
By definition, $[\eta] \in \Ext_\La^1(x'',x') \cap 
\Ima(\beta_{\ii,\cc'',\cc',\f_{x''},\f_{x'}})$
if and only if 
\[
([\eta],0,\ldots,0)=
\beta_{\ii,\cc'',\cc',\f_{x''},\f_{x'}}([\eta_0],\ldots,[\eta_{m-2}]),
\]
in other words if and only if there exists  $([\eta_0],\ldots,[\eta_{m-2}])$  
such that $\iota_{x',1}\circ[\eta_0]=[\eta]$ and
$[\eta_{j-1}]\circ\iota_{x'',j} = \iota_{x',j+1}\circ[\eta_j]$ for all 
$1 \le j \le m-2$.

On the other hand, by Lemma~\ref{imagealpha}
$(\f_{x''},\f_{x'}) \in L^1_{\ii,\cc'',\cc',\eta}$ if and only if for 
$0\leq i\leq m-1$
there exist $[\eta^-_i]\in\Ext^1_\La(x''_i,x'_i)$ such that $[\eta^-_0]=[\eta]$
and $[\eta^-_{i-1}]\circ\iota_{x'',i}=\iota_{x',i}\circ[\eta^-_i]$ for
$1\leq i\leq m-1$. 

Taking into account that for each $i$ either $\iota_{x',i}$ or $\iota_{x'',i}$
is an identity, it is now easy to see that the  two conditions are equivalent.
\end{proof}

%%%%%%%%%%%%%%%%%%%%%%%%%
\subsection{}\label{par3}
%%%%%%%%%%%%%%%%%%%%%%%%%
As before, let $x' \in \La_{\V'}$, $x'' \in \La_{\V''}$, and let
\[
\f_{x'} = 
(x' = x_0' \supseteq x_1' \supseteq  \cdots \supseteq 
x_m' = 0)
\text{\;\;\; and \;\;\;}
\f_{x''} = 
(x'' = x_0'' \supseteq x_1'' \supseteq \cdots \supseteq 
x_m'' = 0)
\]
be flags 
with $(\f_{x'},\f_{x''}) \in \Phi_{\ii,\cc',x'} \times \Phi_{\ii,\cc'',x''}$
for some $(\cc',\cc'') \in W_{\V',\V''}$.

Due to the Calabi-Yau property of $\La$, see~\ref{ssecCY}, 
we have non-degenerate natural pairings
\begin{align*}
\phi_V\df V\times V'\ra\C,\quad 
\left(([\eta_0],\ldots,[\eta_{m-2}]),([\delta_0],\ldots,[\delta_{m-2}])\right)&
\mapsto \sum_{j=0}^{m-2} \phi_{x''_j,x'_{j+1}}([\eta_j])([\delta_j])\\
\phi_W\df W\times W'\ra\C,\quad
\left(([\psi_0],\ldots,[\psi_{m-2}]),([\eps_0],\ldots,[\eps_{m-2}])\right)&
\mapsto \sum_{j=0}^{m-2} \phi_{x''_j,x'_j}([\psi_j])([\eps_j])
\end{align*}
for the spaces $(V,V')$ and $(W,W')$ defined in~\eqref{eqn:betp} 
and~\eqref{eqn:bet}.

\begin{Lem}\label{dualofbeta}
For all $([\eta_0],\ldots,[\eta_{m-2}]) \in V$ and 
$([\eps_0],\ldots,[\eps_{m-2}]) \in W'$ we have
\[
 \phi_V(([\eta_0],\ldots,[\eta_{m-2}]),\beta'(([\eps_0],\ldots,[\eps_{m-2}])))=
\phi_W(\beta(([\eta_0],\ldots,[\eta_{m-2}])),([\eps_0],\ldots,[\eps_{m-2}])).
\]
\end{Lem}

\begin{proof}
By definition, the left hand side of the equation is
\[
\sum_{j=0}^{m-3} \phi_{x''_j,x'_{j+1}}
([\eta_j])([\eps_j]\circ\iota_{x',j+1}-\iota_{x'',j+1}\circ [\eps_{j+1}]) +
\phi_{x''_{m-2},x'_{m-1}}([\eta_{m-1}])([\eps_{m-2}]\circ\iota_{x''_{m-1}}),\\
\]
using~\eqref{eqn:CY-Func}, this is equal to
\[
\phi_{x''_0,x'_0}(\iota_{x',1}\circ[\eta_0])([\eps_0]) + 
\sum_{j=1}^{m-2} \phi_{x''_j,x'_j}
(\iota_{x',j+1}\circ [\eta_j]-[\eta_{j-1}]\circ\iota_{x'',j})([\eps_j]),
\]
which is by definition the right hand side of our claim.
\end{proof}

\begin{Prop}\label{orthogcor3}
For $[\eps] \in \Ext_\La^1(x',x'')$ the
following are equivalent:
\begin{itemize}

\item[(i)]
$[\eps] \in p_0(\Ker(\beta'_{\ii,\cc',\cc'',\f_{x'},\f_{x''}}))$;

\item[(ii)]
$[\eps]$ is orthogonal to $\Ext_\La^1(x'',x') \cap 
\Ima(\beta_{\ii,\cc'',\cc',\f_{x''},\f_{x'}})$.

\end{itemize}
\end{Prop}

\begin{proof}
By Lemma \ref{orthogcor2} and Lemma \ref{dualofbeta} we have
\[
\Ker(\beta'_{\ii,\cc',\cc'',\f_{x'},\f_{x''}})
=
\left(\Ima(\beta_{\ii,\cc'',\cc',\f_{x''},\f_{x'}})\right)^\perp.
\]
Putting $W'_0=\Ext^1_\La(x_0',x_0'')=\Ext_\La^1(x',x'')$
and $W_0=\Ext_\La^1(x_0'',x_0')=\Ext_\La^1(x'',x')$,
it follows that $p_0(\Ker(\beta'_{\ii,\cc',\cc'',\f_{x'},\f_{x''}}))$
is the orthogonal in $W'_0$ of 
$W_0\cap\Ima(\beta_{\ii,\cc'',\cc',\f_{x''},\f_{x'}})$.
\end{proof}

\begin{Cor}\label{orthogcor}
Assume that $\dim \Ext_\La^1(x',x'') = 1$, and let
\[
\eps: \; 0 \to x'' \to x \to x' \to 0 
\text{\;\;\; and \;\;\;}
\eta: \; 0 \to x' \to y \to x'' \to 0
\]
be non-split short exact sequences.
Then the following are equivalent:
\begin{itemize}

\item[(i)]
$(\f_{x'},\f_{x''}) \in L^1_{\ii,\cc',\cc'',\eps}$;

\item[(ii)]
$(\f_{x''},\f_{x'}) \in L^2_{\ii,\cc'',\cc',\eta}$.

\end{itemize}
\end{Cor}

\begin{proof}
By Lemma \ref{betadash} we know that $(i)$ holds if and only if
$[\eps] \in p_0(\Ker(\beta'_{\ii,\cc',\cc'',\f_{x'},\f_{x''}}))$.
By Proposition \ref{orthogcor3} this is equivalent to
$[\eps]$ being orthogonal to $\Ext_\La^1(x'',x') \cap 
\Ima(\beta_{\ii,\cc'',\cc',\f_{x''},\f_{x'}})$.
Since $\dim \Ext_\La^1(x'',x') = 1$ and since $\eps$ is non-split, 
this is the case if and only if
\[
\Ext_\La^1(x'',x') \cap \Ima(\beta_{\ii,\cc'',\cc',\f_{x''},\f_{x'}}) = 0.
\]
By Lemma \ref{beta} and the fact that $\eta$ is non-split this happens 
if and only if
$(\f_{x''},\f_{x'}) \in L^2_{\ii,\cc'',\cc',\eta}$.
\end{proof}

Let us stress that Proposition \ref{orthogcor3} is very similar to
\cite[Proposition 4]{CK}.
We just had to adapt Caldero and Keller's result to our setting 
of preprojective algebras.
Our proof then follows the main line of their proof. 
In fact our
situation turns out to be easier to handle, because we work with
modules and not with objects in cluster categories.

%%%%%%%%%%%%%%%%%%%%%%%%%%%%%%%%%%%%%%%%%%%%%%%%%%%%%%%%%%%%%%%%%%%%%

\section{Euler characteristics of flag varieties}\label{lulemma}

%%%%%%%%%%%%%%%%%%%%%%%%%%%%%%%%%%%%%%%%%%%%%%%%%%%%%%%%%%%%%%%%%%%%%
This section is inspired by the proof of \cite[Lemma 4.4]{Lu}.

%%%%%%%%%%%%%%%%
\subsection{}
Let $\V$, $\W$ and $\T$ be $I$-graded vector spaces with
$\V = \T \oplus \W$. 
We shall identify $\V/\W$ with $\T$.
We write $m=\dim \V$.
Let $\Phi(\V)$ denote the 
variety of complete $I$-graded flags in $\V$.
Let
\[
\pi\df \Phi(\V) \to \Phi_m(\T) \times \Phi_m(\W)
\]
be the morphism which maps a flag $(\V^l)_{0 \le l \le m}$ to 
the pair of flags 
$((\T^l)_{0 \le l \le m}, (\W^l)_{0 \le l \le m})$, 
defined by
$\W^l = \W \cap \V^l$ and $\T^l = p_1(\V^l)$. 
Here 
$$
p_1\df \V = \T \oplus \W \to \T
$$ 
is just the projection onto
$\T$,
and $\Phi_m(\T)$ and $\Phi_m(\W)$ denote the varieties of $m$-step
$I$-graded flags in $\W$ and $\T$, respectively.
Note that
$\T^j/\T^{j+1}$ and  $\W^j/\W^{j+1}$ have dimension 0 or 1 for all $j$.
We want to study the fibers of the morphism $\pi$.

Let $(\f' = (\T^l), \f''=(\W^l)) \in \Ima(\pi)$.
For an $I$-graded linear map $z\df \T \to \W$ we define a flag
\[
\f_z = \f_{z,\f',\f''} = (\V_z^l)_{0 \le l \le m}\in\Phi(\V) 
\]
by 
\[
\V_z^l =  (0,\W^l) + \{ (t,z(t)) \mid t \in \T^l \} \subseteq \T^l\oplus \W.
\]
Clearly, we have $\f_z \in \pi^{-1}(\f',\f'')$.
 
We say that two $I$-graded linear maps $f,g\df \T \to \W$ are
{\it equivalent} if $\V_f^l = \V_g^l$ for all $l$.
In this case we write $f \sim_{\f',\f''} g$ or just $f \sim g$.
In other words, $f$ and $g$ are equivalent if and only
if $(f-g)(\T^l) \subseteq \W^l$ for all $l$. 

The next lemma shows that every flag in $\pi^{-1}(\f',\f'')$
is of the form $\f_z$ for some $z$.

\begin{Lem}
If $\pi(\f) = (\f',\f'')$ then there exists an $I$-graded linear
map $z\df \T \to \W$ such that $\f = \f_{z,\f',\f''}$.
\end{Lem}

\begin{proof}
Choose a decomposition for the two graded vector spaces
$$
\T = \bigoplus_{i=0}^m \T(i) \text{ and } \W = \bigoplus_{i=0}^m \W(i)
$$
such that
$$
\T^l = \bigoplus_{i=l}^m \T(i) \text{ and } \W^l = \bigoplus_{i=l}^m \W(i)
$$
for all $l$.
Put 
$$
\W_c^l := \bigoplus_{i=0}^{l-1}\W(i),
$$
and take a flag $\f = (\V^l)_{0 \le l \le m}$ such that 
$\pi(\f) = (\f',\f'')$.
Since we have a short exact sequence 
\[
0 \to \W^l \to \V^l \to \T^l \to 0,
\]
there exists a unique $I$-graded linear map
$w^l\df \T^l \ra \W_c^l$ such that
\[ 
\V^l = (0,\W^l) + \{ (t,w^l(t)) \mid t \in \T^l \}. 
\] 
The conditions $\V^l\supset\V^{l+1}$ and $\W_c^l \subseteq \W_c^{l+1}$ imply 
\begin{equation} \label{eq:lem321}
w^{l+1}(t)-w^l(t) \in \W_c^{l+1} \cap \W^l = \W(l)
\end{equation}
for all $t\in\T^{l+1}$.
Hence, for $t\in\T^{l+1}$ we have
\[
(t,w^{l+1}(t)) - (t,w^{l}(t))= (0,w^{l+1}(t)-w^l(t)) \in \W^l \subseteq \V^l.
\]
Now, define $z\df\T\ra\W$ on the summands of $\T$ by
$z(t)=w^{i}(t)$ for $t\in\T(i)$. 
\end{proof}

Thus the fiber
$\pi^{-1}(\f',\f'')$
can be identified with the vector space 
\[
\Hom(\T,\W)/\!\sim_{\f',\f''}
\] 
of $I$-graded linear maps $\T \to \W$ up to equivalence.

%%%%%%%%%%%%%%%%
\subsection{}
Let $\Ga = \C Q'/J$ be an algebra where $Q'$ is a finite quiver 
with set of vertices $I$ 
and $J$ is an ideal in the path algebra $\C Q'$.
For a finite-dimensional $I$-graded vector space $\U$ let 
$\Ga_\U$ be the affine variety of $\Ga$-module structures on $\U$.
Fix a short exact sequence
\[
\eps: \;\;\; 0 \to x'' \to x \to x' \to 0
\]
of $\Ga$-modules with
$x \in \Ga_\V$, $x' \in \Ga_\T$ and $x'' \in \Ga_\W$.
To an $I$-graded linear map $z\df \T \to \W$ we associate linear maps
$z^l\df \T^l \to \W/\W^l$ defined by 
$z^l(t) = z(t) + \W^l$.
For a vertex $i \in I$ let $z_i\df \T_i \to \W_i$ be the degree $i$ part of
$z$. 
The next lemma follows from the definitions.

\begin{Lem}\label{eulerlemma1}
If $x = x' \oplus x''$ and if $\f'$ and $\f''$ are flags of submodules
of $x'$ and $x''$, respectively, then 
$\f_{z,\f',\f''}$ is a flag of submodules of
$x' \oplus x''$ if and only if the map $z^l$ is a 
module homomorphism for all $l$.
\end{Lem}

We now deal with the case when the short exact sequence $\eps$
does not split.
A module $m \in \Ga_\V$ can be interpreted as a tuple 
$m = (m(a))_a$ where for each arrow $a \in Q'_1$ we have a linear map
$m(a)\df \V_{s(a)} \to \V_{e(a)}$ where
$s(a)$ and $e(a)$ denote the start and end vertex of the arrow $a$.
Thus given our short exact sequence
\[
\eps: \;\;\; 0 \to x'' \to x \to x' \to 0
\]
we can assume without loss of generality that 
the linear maps $x(a)$ are of the form
\[
x(a) = \left( \begin{matrix} x'(a) & 0\\ y(a) & x''(a) \end{matrix} \right)
\]
where $y(a)\df \T_{s(a)} \to \W_{e(a)}$ are certain linear maps.
The proof of the following statement is again straightforward, compare also
the proof of \cite[Lemma 4.4]{Lu}.

\begin{Lem}\label{eulerlemma2}
If $\f'$ and $\f''$ are flags of submodules
of $x'$ and $x''$, respectively, then 
$\f_{z,\f',\f''}$ is a flag of submodules of
$x$ if and only if 
\[
(x''(a) z_{s(a)} - z_{e(a)} x'(a) - y(a))(\T^l_{s(a)}) \subseteq \W^l_{e(a)}
\]
for all $0 \le l \le m-1$ and $a \in Q_1'$.
\end{Lem}

%%%%%%%%%%%%%
\subsection{}
%%%%%%%%%%%%%
Now we apply the above results to the case of the preprojective
algebra $\La$.

\begin{Lem}\label{fibredim}
For a short exact sequence $\eps: \; 0 \to x'' \to x \to x' \to 0$ and 
$(\cc',\cc'') \in W_{\V',\V'',\eps}$
the fibers of the morphism
\[
\alpha_{\ii,\eps}(\cc',\cc'')\df
\Phi_{\ii,x}(\cc',\cc'',\eps)
\to \Phi_{\ii,\cc',x'}\times \Phi_{\ii,\cc'',x''} 
\]
are isomorphic to affine spaces, moreover 
$\Ima(\alpha_{\ii,\eps}(\cc',\cc''))=L^1_{\ii,\cc',\cc'',\eps}$.
\end{Lem}

\begin{proof}
Both conditions in Lemma \ref{eulerlemma1} and Lemma \ref{eulerlemma2}
are linear. The last equation follows from the definitions, 
see~\ref{ssec:refine}.
\end{proof}

\begin{Cor}\label{eulerformula2}
For  a short exact sequence $\eps: \; 0 \to x'' \to x \to x' \to 0$ 
and $(\cc',\cc'') \in W_{\V',\V''}$
we have
\[
\chi(\Phi_{\ii,x}(\cc',\cc'',\eps)) = \chi(L^1_{\ii,\cc',\cc'',\eps}).
\] 
\end{Cor}

\begin{proof}
This follows from Lemma \ref{fibredim} and Proposition \ref{eulerfibres}.
\end{proof}

\begin{Cor}\label{eulerformula1}
If $x = x' \oplus x''$, then for all short exact sequences
$\eps: \; 0 \to x'' \to x \to x' \to 0$ and 
$(\cc',\cc'') \in W_{\V',\V''}$ we have
\[
\chi(\Phi_{\ii,x}(\cc',\cc'',\eps))=
\chi(\Phi_{\ii,\cc',x'}\times\Phi_{\ii,\cc'',x''})
=\chi( \Phi_{\ii,\cc',x'})\cdot \chi(\Phi_{\ii,\cc'',x''}).
\]
\end{Cor}

\begin{proof}
If $x = x' \oplus x''$ and $(\cc',\cc'') \in W_{\V',\V''}$, then 
$L_{\ii,\cc',\cc'',\eps}^1 = \Phi_{\ii,\cc',x'} \times 
\Phi_{\ii,\cc'',x''}$ for all $\eps$, see \cite{GLS}.
\end{proof}

In \cite{GLS} we claimed that for $(\cc',\cc'') \in W_{\V',\V'',\eps}$ 
the map 
$\Phi_{\ii,x' \oplus x''}(\cc',\cc'',\eps) \to 
\Phi_{\ii,\cc',x'} \times \Phi_{\ii,\cc'',x''}$
is a vector bundle, referring to \cite[Lemma 4.4]{Lu}.
What we had in mind was an argument as above.
However this only proves that the fibers of this map are affine spaces.
Nevertheless, for the Euler characteristic calculation needed in \cite{GLS}
this is enough because of Proposition \ref{eulerfibres}.

%%%%%%%%%%%%%%%%%%%%%%%%%%%%%%%%%%%%%%%%%%%

\section{Proof of Theorem \ref{mainresult}}

%%%%%%%%%%%%%%%%%%%%%%%%%%%%%%%%%%%%%%%%%%%

Assume $\dim \Ext_\La^1(x',x'') = 1$, and let
\[
\eps: \; 0 \to x'' \to x \to x' \to 0 
\text{\;\;\; and \;\;\;}
\eta: \; 0 \to x' \to y \to x'' \to 0
\]
be non-split short exact sequences.
We obtain the following diagram of maps:
\[
\xymatrix{
\Phi_{\ii,x}(\cc',\cc'',\eps) 
\ar[d]_{\alpha_{\ii,\eps}(\cc',\cc'')} & & 
\Phi_{\ii,y}(\cc'',\cc',\eta) 
\ar[d]^{\alpha_{\ii,\eta}(\cc'',\cc')}\\
L^1_{\ii,\cc',\cc'',\eps} 
\ar[d]_{\iota_{\ii,\cc',\cc'',\eps}} && 
L^1_{\ii,\cc'',\cc',\eta}  \ar[d]^{\iota_{\ii,\cc'',\cc',\eta}}\\
L^1_{\ii,\cc',\cc'',\eps} \cup L^2_{\ii,\cc',\cc'',\eps} \ar@{=}[d] 
&& L^1_{\ii,\cc'',\cc',\eta} \cup L^2_{\ii,\cc'',\cc',\eta} 
\ar@{=}[d]\\
\Phi_{\ii,\cc',x'} \times \Phi_{\ii,\cc'',x''} \ar[rr]^i 
& &
\Phi_{\ii,\cc'',x''} \times \Phi_{\ii,\cc',x'}
} 
\]
Here, the map $i$ is the isomorphism which maps $(\f_{x'},\f_{x''})$ to
$(\f_{x''},\f_{x'})$, and 
$\iota_{\ii,\cc',\cc'',\eps}$ and
$\iota_{\ii,\cc'',\cc',\eta}$ are the natural inclusion maps.
By Corollary \ref{eulerformula2}, 
\[
\chi(\Phi_{\ii,x}(\cc',\cc'',\eps))
= 
\chi(L^1_{\ii,\cc',\cc'',\eps})
\text{\;\;\; and \;\;\;}
\chi(\Phi_{\ii,y}(\cc'',\cc',\eta))
= 
\chi(L^1_{\ii,\cc'',\cc',\eta}).
\] 
We know from Corollary \ref{eulerformula1} that
\[
\chi(\Phi_{\ii,x' \oplus x''}(\cc',\cc'',\theta)) = 
\chi(\Phi_{\ii,\cc',x'} \times \Phi_{\ii,\cc'',x''}) = 
\chi(L^1_{\ii,\cc',\cc'',\eps}) + \chi(L^2_{\ii,\cc',\cc'',\eps}), 
\]
where $\theta$ is any (split) short exact sequence of the form
$0 \to x'' \to x'' \oplus x' \to x' \to 0$.
By Corollary \ref{orthogcor} the isomorphism $i$ induces isomorphisms
\[
i_{1,2}\df L_{\ii,\cc',\cc'',\eps}^1 \to L_{\ii,\cc'',\cc',\eta}^2
\text{\;\;\; and \;\;\;}
i_{2,1}\df L_{\ii,\cc',\cc'',\eps}^2 \to L_{\ii,\cc'',\cc',\eta}^1
\]
which implies
$\chi(L^2_{\ii,\cc',\cc'',\eps}) = \chi(L^1_{\ii,\cc'',\cc',\eta})$.
Combining these facts we get
\begin{align*}
\chi(\Phi_{\ii,x' \oplus x''}) 
&= \sum_{(\cc',\cc'') \in W_{\V',\V''}} 
\chi(\Phi_{\ii,x' \oplus x''}(\cc',\cc'',\theta))\\
&= \sum_{(\cc',\cc'') \in W_{\V',\V''}}
\chi(L^1_{\ii,\cc',\cc'',\eps}) +
\sum_{(\cc',\cc'') \in W_{\V',\V''}} 
\chi(L^2_{\ii,\cc',\cc'',\eps})\\
&= \sum_{(\cc',\cc'') \in W_{\V',\V''}} 
\chi(L^1_{\ii,\cc',\cc'',\eps}) +
\sum_{(\cc',\cc'') \in W_{\V',\V''}} 
\chi(L^1_{\ii,\cc'',\cc',\eta})\\
&= \sum_{(\cc',\cc'') \in W_{\V',\V''}} 
\chi(\Phi_{\ii,x}(\cc',\cc'',\eps)) +
 \sum_{(\cc',\cc'') \in W_{\V',\V''}} 
\chi(\Phi_{\ii,y}(\cc'',\cc',\eta))\\
&= \chi(\Phi_{\ii,x}) + \chi(\Phi_{\ii,y}).
\end{align*}
By the considerations in Section \ref{reformulation} 
this finishes the proof of
Theorem \ref{mainresult}.

%%%%%%%%%%%%%%%%%%%%%%%%%%%%%%%%%%%%%%%%%%%%%%%%%%%%%%%%%%%%%%

\section{The general case}\label{generalcase}

%%%%%%%%%%%%%%%%%%%%%%%%%%%%%%%%%%%%%%%%%%%%%%%%%%%%%%%%%%%%%%

%%%%%%%%%%%%%%%%%%%%%%%%
\subsection{Derivations}
%%%%%%%%%%%%%%%%%%%%%%%%
For $\La$-modules $x'$ and $x''$ let
$D_\La(x',x'')$ be the vector space of all tuples 
$d = (d(b))_{b \in \overline{Q}_1} $ of linear maps
$d(b) \in \Hom_\C(x'_{s(b)},x''_{e(b)})$ such that
\begin{equation}\label{derivationeq}
\sum_{a \in Q_1: s(a)=p} (d(\bar{a})x'(a) + x''(\bar{a})d(a)) -
\sum_{a \in Q_1: e(a)=p} (d(a)x'(\bar{a}) + x''(a)d(\bar{a})) = 0
\end{equation}
for all $p \in Q_0$.
We call the elements in $D_\La(x',x'')$ {\it derivations}.

Let $d = (d(a))_{a \in \overline{Q}_1}$ with
$d(a) \in \Hom_\C(x'_{s(a)},x''_{e(a)})$, and for
each $a \in \overline{Q}_1$ let 
\[
E_d(a) = \left( \begin{matrix}
x'(a) & 0\\
d(a) & x''(a)
\end{matrix} \right).
\]
Then $E_d = (E_d(a))_{ a \in \overline{Q}_1}$ defines a $\La$-module if and
only if the maps $d(a)$ satisfy Equation (\ref{derivationeq}) for all $p$.
In this case we obtain an obvious short exact sequence
\[
\eps_d: \;\;\; 0 \to x'' \to E_d \to x' \to 0.
\]

%%%%%%%%%%%%%%%%%%%%%%%%%%%%%%
\subsection{Inner derivations}
%%%%%%%%%%%%%%%%%%%%%%%%%%%%%%
For $\La$-modules $x'$ and $x''$ let
$I_\La(x',x'')$ be the vector space of all tuples 
$i = (i(b))_{b \in \overline{Q}_1} $ of linear maps
$i(b) \in \Hom_\C(x'_{s(b)},x''_{e(b)})$ such that
for some $(\phi_q)_{q \in Q_0}$ with 
$\phi_q \in \Hom_\C(x'_q,x''_q)$ we have
\[
i(b) = \phi_{e(b)}x'(b) - x''(b)\phi_{s(b)}
\]
for all $b \in \overline{Q}_1$.
The elements in $I_\La(x',x'')$ are called {\it inner derivations}
and we have obviously $I_\La(x',x'') \subseteq D_\La(x',x'')$.
Let 
$$
\pi\colon D_\La(x',x'')\to \Ext_\La^1(x',x'')
$$ 
be defined by
$\pi(d) = [\eps_d]$. 
It is known that the kernel of $\pi$ is just the set $I_\La(x',x'')$ of inner
derivations, hence we obtain an exact sequence 
\[
0 \to \Hom_\La(x',x'') \to \bigoplus_{q \in Q_0} \Hom_\C(x'_q,x''_q) 
\to D_\La(x',x'') \xrightarrow{\pi} \Ext_\La^1(x',x'') \to 0.
\]
We choose a fixed vector space decomposition
\[
D_\La(x',x'') = I_\La(x',x'') \oplus E_\La(x',x'').
\]
We can therefore identify $\Ext_\La^1(x',x'')$ with $E_\La(x',x'')$.
We also can identify $\PP \Ext_\La^1(x',x'')$ with $\PP E_\La(x',x'')$.
Set 
$$
E^*_\La(x',x'') = E_\La(x',x'') \setminus \{0 \}.
$$

%%%%%%%%%%%%%%%%%%%%%%%%%%%%%%%%%%%%%%%%%%%%%%%%%%%%%%%%
\subsection{Principal $\C^*$-bundles}\label{ssec:pribdl}
%%%%%%%%%%%%%%%%%%%%%%%%%%%%%%%%%%%%%%%%%%%%%%%%%%%%%%%%
Recall that the multiplicative group $\C^*$ is {\em special} in the sense that
each principal $\C^*$-bundle is locally trivial in the Zariski topology,
see \cite[\S 4]{Serre}. 
Thus, if $\pi\df P\ra Q$ is such a bundle, $\pi$ admits local sections. 
Since moreover the action of $\C^*$ on $P$ is free by
definition, we conclude that $(\pi, Q)$ is a geometric quotient for the
action of $\C^*$ on $P$, see for example \cite[Lemma 5.6]{Bonga}. 
Note moreover that $\pi$ is flat (and in particular open) since  
locally it is just a projection.

As a rule, we will write $\C^* x$ for the $\C^*$-orbit of $x$ if $x$ belongs
to a principal $\C^*$-bundle.

%%%%%%%%%%%%%%%%%%%%%%%%%%%%%%%%%%%%%%%%%%%%%
\subsection{The varieties ${EF}_\ii(x',x'')$}
%%%%%%%%%%%%%%%%%%%%%%%%%%%%%%%%%%%%%%%%%%%%%
Let $x'\in\La_{\V'}$ and $x''\in\La_{\V''}$ and
$\ii=(i_1,i_2,\ldots, i_m)$ be a word in $I^m$ such that 
$\abs{\ii}=\dimv (x'\oplus x'')$.
The action of $\C^*$ on $\V'\oplus\V''$ defined by 
$$
\lam\star(v',v''):= (v',\lam\cdot v'')
$$
induces an action on the flag variety $\Phi_\ii(\V'\oplus\V'')$:
if 
$$
x^{\bullet}=(x^0\supseteq x^1\supseteq \cdots \supseteq x^m),
$$ 
then the $i$-th component of
$(\lam\star x^{\bullet})$ is $\{\lam\star z\mid z\in x^i\}$. 

On the other hand we have the principal $\C^*$-bundle 
$E^*_\La(x',x'')\ra \PP E_\La(x',x'')$. 
Thus we obtain by \cite[\S 3.2]{Serre} a new principal $\C^*$-bundle
\[
\tilde{q}\df E^*_\La(x',x'')\times \Phi_\ii(\V'\oplus \V'')\ra 
E^*_\La(x',x'')\times^{\C^*} \Phi_\ii(\V'\oplus \V'').
\]
We consider 
\[
\widetilde{EF}_\ii(x',x''):= 
\{(d,x^{\bullet})\in E^*_\La(x',x'')\times \Phi_\ii(\V'\oplus \V'')\mid
x^{\bullet} \in \Phi_{\ii, E_d} \}.
\]
This is clearly a closed subset of 
$E^*_\La(x',x'')\times\Phi_\ii(\V'\oplus \V'')$,
and it is moreover $\C^*$-stable since
\begin{equation}\label{eq:iso}
\left( \begin{matrix} \id_{\V'} & 0\\ 0 & \lam \id_{\V''} \end{matrix} 
\right)\colon 
E_d = \left( \begin{matrix} x'&0\\d&x'' \end{matrix} \right) \to
E_{\lam d} = \left( \begin{matrix} x'&0\\\lam d & x'' \end{matrix} \right)
\end{equation}
is an isomorphism of $\La$-modules.
We conclude that 
$$
EF_{\ii}(x',x''):=\tilde{q}(\widetilde{EF}_\ii(x',x''))
$$ 
is closed
in $E^*_\La(x',x'')\times^{\C^*} \Phi_\ii(\V'\oplus \V'')$ since $\tilde{q}$
is open, see \ref{ssec:pribdl}. 
Thus, the restriction of $\tilde{q}$ to
\[
q\df \widetilde{EF}_\ii(x',x'')\ra EF_{\ii}(x',x'')
\]
is again a principal $\C^*$-bundle \cite[\S 3.1]{Serre}, and in particular
$(q,EF_\ii(x',x''))$ is a geometric quotient for the action of $\C^*$ on
$\widetilde{EF}_\ii(x',x'')$.

The projection 
\[
\widetilde{p}_1\df\widetilde{EF}_{\ii}(x',x'') \ra \PP E_\La(x',x''), 
(d,x^{\bullet})\mapsto \C^* d
\]
is constant on $\C^*$-orbits. So we 
obtain a morphism $p_1\df EF_\ii(x',x'') \to {\PP}E_\La(x',x'')$
which maps $\C^*(d,x^{\bullet})$ to $\C^* d$.
In other words, we have 
$$
p_1^{-1}(\C^*d) \cong \Phi_{\ii,E_d}.
$$

%%%%%%%%%%%%%%%%%%%%%%%%%%%%%%%%%%%%%%%%%%%%%%%%%%%%%%%%%%%%%%%%%%%%%
\subsection{$\delta$-stratification of $EF_\ii(x',x'')$}\label{E-del}
%%%%%%%%%%%%%%%%%%%%%%%%%%%%%%%%%%%%%%%%%%%%%%%%%%%%%%%%%%%%%%%%%%%%%
Let $\ee := \dimv(x') +\dimv(x'')$ and consider the 
$\delta$-stra\-tifi\-cation 
\[
\La_\ee = \bigsqcup_{x\in R(\ee)} \str{x}. 
\]
We claim that 
\[
EF_\ii(x',x'')_{\str{x}} := \{ \C^*(d,x^{\bullet}) \in {EF}_\ii(x',x'') 
\mid E_d \in \str{x} \}
\]
is a constructible set. Indeed, 
we have a morphism $\iota\df E_\La(x',x'') \to \La_\ee$ which maps
$d$ to $E_d$. 
Since $\str{x}$ is constructible and $\Gl_\ee$-invariant, 
$$
E_\La(x',x'')_{\str{x}} = \{ d \in E_\La(x',x'') \mid E_d \in \str{x}\}
$$ 
is constructible
and $\C^*$-invariant (see also \eqref{eq:iso}). 
Now,
$$
EF_\ii(x',x'')_{\str{x}} = p_1^{-1}(\PP E_\La(x',x'')_{\str{x}})
$$ 
is also constructible.

The fibers $p_1^{-1}([d])$ are identified with the varieties $\Phi_{\ii,E_d}$. 
Since they all
come from the same $\delta$-stratum, they all have the same Euler 
characteristic.
Thus from Proposition \ref{eulerfibres} we obtain 
\[
\chi({EF}_\ii(x',x'')_{\str{x}}) 
= \chi(\PP E_\La(x',x'')_{\str{x}})
\cdot \chi(\Phi_{\ii,E_d})
= \chi(\PP\Ext^1_\La(x',x'')_{\str{x}})
\cdot \chi(\Phi_{\ii,E_d})
\]
for any $d \in E_\La(x',x'')_{\str{x}}$. 
It follows that,
\begin{equation}\label{eq:main1}
\chi(EF_\ii(x',x''))=
\sum_{x\in R(\ee)} \chi(\PP\Ext^1_\La(x',x'')_{\str{x}})\cdot \chi(\Phi_{\ii,x}).
\end{equation}

%%%%%%%%%%%%%%%%%%%%%%%%%%%%%%%%%%%%%%%%
\subsection{Proof of Theorem \ref{thm1}}
%%%%%%%%%%%%%%%%%%%%%%%%%%%%%%%%%%%%%%%%
We define a morphism of varieties
\[
\pi\df {EF}_\ii(x',x'') \to \bigsqcup_{(\cc',\cc'')\in W_{\V',\V''}} 
(\PP E_\La(x',x'') \times \Phi_{\ii,\cc',x'} \times \Phi_{\ii,\cc'',x''})
\]
which maps $\C^*(d,x^{\bullet})$ to $(\C^* d,\f_{x'},\f_{x''})$
where $(\f_{x'},\f_{x''})$ is the pair of flags in $x'$ and $x''$ induced by
$x^{\bullet}$ via the short exact sequence
\[
\eps_d: \;\;\; 0 \to x'' \to E_d \to x' \to 0.
\]
(Observe that $(d,x^{\bullet})$ and $(\lam d,\lam \star x^{\bullet})$ induce 
the
same pair of flags via the sequences $\eps_d$ and $\eps_{\lam d}$,
respectively.)
We denote by $L^1$ the image of $\pi$, and by $L^2$ the complement
of $L^1$ in 
$
\bigsqcup_{(\cc',\cc'') \in W_{\V',\V''}} 
(\PP E_\La(x',x'') \times \Phi_{\ii,\cc',x'} \times \Phi_{\ii,\cc'',x''}).
$
The following diagram illustrates the situation:

$$
\unitlength1.0cm
\begin{picture}(14,4)
\put(4,4){$\widetilde{EF}_\ii(x',x'') \subseteq
E_\La^*(x',x'') \times \Phi_\ii(\V' \oplus \V'')$}
\put(4.0,3.8){\vector(-1,-1){1.4}}
\put(4.6,3.8){\vector(0,-1){1.4}}
\put(4.8,3.1){$q$}
\put(2.8,3.1){$\widetilde{p_1}$}
\put(8.3,3.8){\vector(0,-1){1.4}}
\put(8.5,3.1){$\widetilde{q}\colon (d,x^\bullet) \mapsto \C^*(d,x^\bullet)$}
\put(0.6,2){$\PP E_\La(x',x'')$}
\put(4,2){${EF}_\ii(x',x'') \subseteq
E_\La^*(x',x'') \times^{\C^*} \Phi_\ii(\V' \oplus \V'')$}
\put(3.8,2.1){\vector(-1,0){1.2}}
\put(3.2,1.8){$p_1$}
\put(4.6,1.8){\vector(0,-1){1.4}}
\put(4.8,1.1){$\pi\colon \C^*(d,x^\bullet) \mapsto (\C^* d,\f_{x'},\f_{x''})$}
\put(4,0){$L^1 \sqcup L^2 = \bigsqcup_{(\cc',\cc'') \in W_{\V',\V''}} 
(\PP E_\La(x',x'') \times \Phi_{\ii,\cc',x'} \times \Phi_{\ii,\cc'',x''})$}
\end{picture}
$$

\begin{Prop} 
The following hold:
\begin{itemize}

\item[(a)]
$\chi(L^1)=\chi(EF_{\ii}(x',x''))$;

\item[(b)]
$
\chi\left(\bigsqcup_{(\cc',\cc'')\in W_{\V',\V''}} 
(\PP E_\La(x',x'') \times \Phi_{\ii,\cc',x'} \times
\Phi_{\ii,\cc'',x''})\right) = 
\chi(\PP E_\La(x',x'')) \cdot \chi(\Phi_{\ii,x' \oplus x''});
$

\item[(c)]
$\chi(L^2) = \chi({EF}_\ii(x'',x'))$.
\end{itemize}
\end{Prop}

\begin{proof} 
(a) By the same argument as in  Lemma \ref{fibredim}, the fibers of
$\pi$ are isomorphic to affine spaces.
This implies (a) by Proposition \ref{eulerfibres}.

(b) We have to show that
\[
\sum_{(\cc',\cc'')\in W_{\V',\V''}}
\chi(\Phi_{\ii,\cc',x'}) \cdot \chi(\Phi_{\ii'',\cc'',x''})
=\chi(\Phi_{\ii,x'\oplus x''}).
\]
This is explained in the proof of~\cite[Lemma 6.1]{GLS}, 
see also the remark after Corollary \ref{eulerformula1}.

(c) 
Let $(\PP[d],\f_{x'},\f_{x''})\in L^2$.
By Lemma \ref{betadash}, this means that for all $(\cc',\cc'')\in
W_{\V',\V''}$,
\[
d\not \in p_0(\Ker(\beta'_{\ii,\cc',\cc'',\f_{x'},\f_{x''}})).
\] 
By Proposition \ref{orthogcor3}, this means that for all $(\cc',\cc'')\in
W_{\V',\V''}$, 
\[
d \not\perp (E_\La(x'',x') \cap \Ima(\beta_{\ii,\cc'',\cc',\f_{x''},\f_{x'}})).
\]
Therefore there exists 
$d'\in E_\La(x'',x') \cap
\Ima(\beta_{\ii,\cc'',\cc',\f_{x''},\f_{x'}})$
which is not orthogonal to $d$, and a flag $y^\bullet$ of submodules
of the module
\[
\left( \begin{matrix}x''&0\\d'&x' \end{matrix} \right)
\]
such that $y^\bullet$ induces $\f_{x''}$ and $\f_{x'}$. 
We are thus led to consider the constructible set $C$ of all pairs
\[
((\C^* d,\f_{x'},\f_{x''}),\C^*(d',y^\bullet)) \in 
L^2 \times {EF}_\ii(x'',x')
\]
such that $(d',y^\bullet)$ induces $(\f_{x''},\f_{x'})$ and $d$ and $d'$ 
are not
orthogonal for the pairing between $E_\La(x',x'')$ and
$E_\La(x'',x')$. 
Let us consider the two natural projections:
\[
\xymatrix{
&C \ar[dl]_{pr_1} \ar[dr]^{pr_2}\\
L^2 && {EF}_\ii(x'',x')
}
\]
We are going to show that all fibers of both projections
have Euler characteristic equal to $1$.
More precisely we have: 
\begin{itemize}

\item[(i)]
The map $pr_1$ is surjective with fibers being extensions of affine
spaces;

\item[(ii)]
The map $pr_2$ is surjective with fibers being affine spaces 
(of constant dimension $\dim E_\La(x'',x')-1$).

\end{itemize}
Let us prove (i). 
Let $(\C^*d,\f_{x'},\f_{x''})\in L^2$.
Let $E_\La(x'',x')_{(\f_{x''},\f_{x'})}$ 
be the set of all $d' \in E_\La(x'',x')$
such that there exists a filtration $y^\bullet$ of
the module 
\[
\left( \begin{matrix}x''&0\\d'&x' \end{matrix} \right)
\]
which induces $\f_{x''}$ and $\f_{x'}$.
This is a vector space, and by the above discussion we know that 
it is not contained in the hyperplane $d^\perp$.
Thus 
$$
d^\perp \cap E_\La(x'',x')_{(\f_{x''},\f_{x'})}
$$ 
is a hyperplane
in $E_\La(x'',x')_{(\f_{x''},\f_{x'})}$, and 
\[
Z := \PP(E_\La(x'',x')_{(\f_{x''},\f_{x'})} \setminus d^\perp) 
\]
is an affine space.
We get a map
\[
pr_1^{-1}(\C^* d ,\f_{x'},\f_{x''}) \to Z
\]
which maps $((\C^* d,\f_{x'},\f_{x''}),\C^*(d',y^\bullet))$ to $\C^*d'$.
By construction this map is surjective and
its fibers are affine spaces which implies  (i).

Let us prove (ii).
Let $\C^*(d',y^\bullet)\in EF_{\ii}(x'',x')$.
Let $\f_{x''}$ and $\f_{x'}$ be the flags induced by $y^\bullet$
on $x''$ and $x'$, respectively.
Then by Lemma \ref{beta} and Proposition \ref{orthogcor3}, 
\[
d'\in p_0(\Ker(\beta'_{\ii,\cc',\cc'',\f_{x'},\f_{x''}}))^\perp
\]
for some $(\cc',\cc'')\in W_{\V',\V''}$.
So if $d \not\in {^\perp}d'$, then $(\C^* d,\f_{x'},\f_{x''})\in L^2$,
by Lemma \ref{betadash}.
Therefore $pr_2^{-1}([d',y^\bullet])$  
can be identified with the projectivization of the set of all 
$d \in E_\La(x',x'')$ such that $d \not\in {^\perp}d'$.
But the projectivization of the complement of a hyperplane in a
vector space of dimension $m$ is an affine space of dimension $m-1$.
\end{proof}

Now, the proposition together with Equation \eqref{eq:main1} imply obviously
Equation \eqref{eqn:ref1}, which is equivalent to Theorem \ref{thm1}.

%%%%%%%%%%%%%%%%%%%%%%%%%%%%%%%%%%%%%%%%%%%%%%%%%%%%%%%%%%%%%

\section{An example}\label{examples}

%%%%%%%%%%%%%%%%%%%%%%%%%%%%%%%%%%%%%%%%%%%%%%%%%%%%%%%%%%%%%

In this section, $\La$ will be the preprojective algebra of type $\Dd_4$,
with underlying quiver 
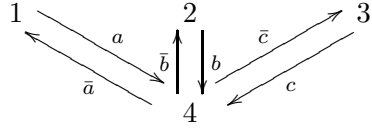
\[
\xymatrix{
1\ar@<1ex>[rd]^a&2\ar@<1ex>[d]^b& 3\ar@<1ex>[ld]^c\\
&{\phantom{XXX}4\phantom{XXX}}
\ar@<1ex>[lu]^{\bar{a}}\ar@<1ex>[u]^{\bar{b}}\ar@<1ex>[ru]^{\bar{c}}
}
\]
We study extensions between 
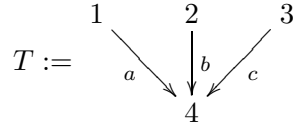
\[
T:=\vcenter{\xymatrix{
1\ar[rd]_a&2\ar[d]^b&3\ar[ld]^c\\
&4
}}
\]
and the simple module $S_4$. We have $\dim\Ext^1_\La(T,S_4)=2$.
The middle terms of the non-split short exact sequences
\[
 0\ra T\ra E\ra S_4\ra 0
\]
are of the form
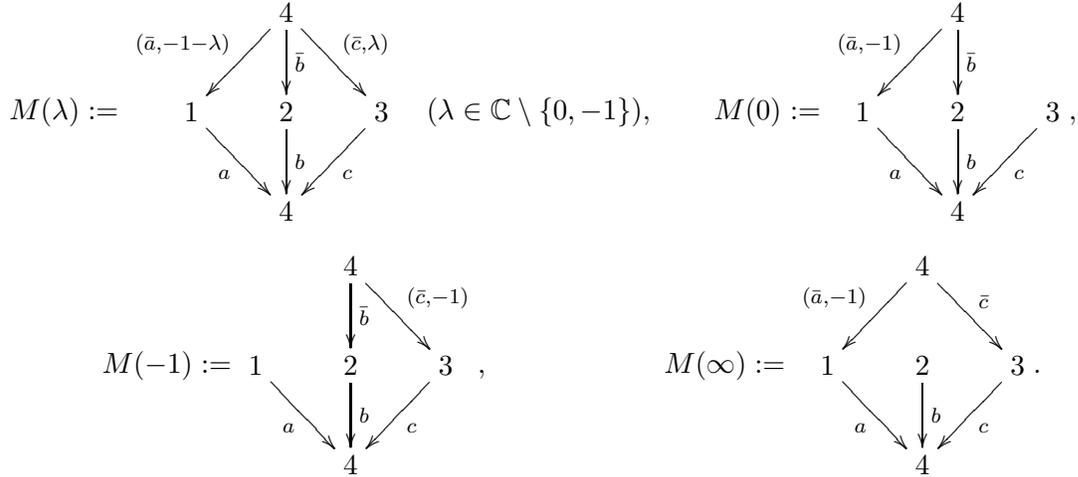
\[
M(\lambda):= 
\vcenter{\xymatrix{
&4\ar[ld]_{(\bar{a},-1-\lambda)}\ar[d]^{\bar{b}}\ar[rd]^{(\bar{c},\lambda)}\\
1\ar[rd]_a&2\ar[d]^b&3\ar[ld]^c\\
&4}}\quad (\lambda\in\C \setminus \{0,-1\}),\qquad
M(0):=
\vcenter{\xymatrix{
&4\ar[ld]_{(\bar{a},-1)}\ar[d]^{\bar{b}}\\
1\ar[rd]_a&2\ar[d]^b&3\ar[ld]^c\\
&4}},\qquad
\]
\[
M(-1):=\vcenter{\xymatrix{
&4\ar[d]^{\bar{b}}\ar[rd]^{(\bar{c},-1)}\\
1\ar[rd]_a&2\ar[d]^b&3\ar[ld]^c\\
&4}},\qquad\qquad\qquad
M(\infty):=\vcenter{\xymatrix{
&4\ar[ld]_{(\bar{a},-1)}\ar[rd]^{\bar{c}}\\
1\ar[rd]_a&2\ar[d]^b&3\ar[ld]^c\\
&4}}.
\]
The middle terms of the non-split short exact sequences
\[
0\ra S_4\ra E\ra T\ra 0
\]
are of the form
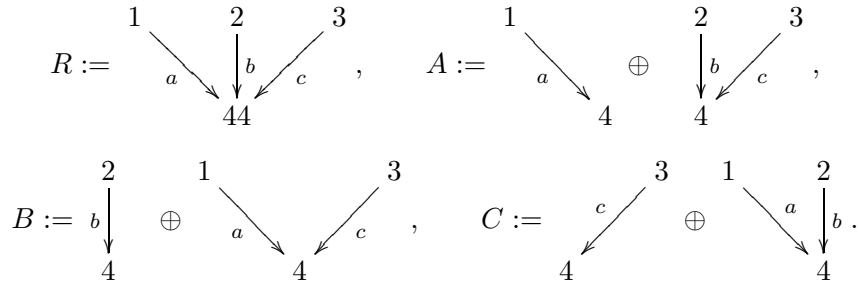
\[
R:= \vcenter{\xymatrix{
1\ar[rd]_a&2\ar[d]^b&3\ar[ld]^c\\
&44
}},\qquad
A:=\vcenter{\xymatrix{1\ar[rd]_a\\&4}}
\oplus\quad\vcenter{\xymatrix{2\ar[d]^b&3 \ar[ld]^c\\4}},
\]
\[
B:=\vcenter{\xymatrix{2\ar[d]_b\\4}}\quad
\oplus\vcenter{\xymatrix{1\ar[rd]_a&&3 \ar[ld]^c\\&4}},\qquad
C:=\vcenter{\xymatrix{&3\ar[ld]_c\\4}}
\oplus\vcenter{\xymatrix{1\ar[rd]^a&2 \ar[d]^b\\&4}}.
\]
Note that $R$ is rigid, and $A,B,C$ belong to the orbit closure of $R$.
Our multiplication formula yields for any $\lam\in\C\setminus\{0,-1\}$
$$
2 \del_T\cdot\del_{S_4} = \left((-1)\del_{M(\lam)}+\del_{M(0)}+
\del_{M(-1)}+\del_{M(\infty)}\right)
+ \left((-1) \del_R +\del_A+ \del_B+\del_C\right).
$$
To see this, note that $\PP\Ext^1_\La(S_4,T)$ is a projective line, with
points identified to the middle terms of short exact sequences
\[
0\ra T \ra M(\lam) \ra S_4\ra 0\quad\text{with } \lam\in\C\cup\{\infty\}.
\]
This becomes ``stratified'' according to the Euler characteristics of
flags of submodules into 
$(\C\setminus\{0,1\})\cup\{0\}\cup\{1\}\cup\{\infty\}$ and we obtain the first
term on the right hand side. Also $\Ext^1_\La(T,S_4)$ is a projective line.
In this situation there are only $4$ isomorphism classes 
$R,A,B,C$ of middle terms for non-split short exact sequences
\[
0\ra S_4\ra X\ra T\ra 0.
\]
However, in the first case we have a $(\C\setminus\{0,1\})$-family of possible
embeddings of $S_4$ into $R$ such that the quotient is $T$ while in the
other three cases there is a unique embedding of this type. This gives the
second term on the right hand side of our equation.

Using the rigid modules
\[
F:=\vcenter{\xymatrix{
1\ar[rd]_a\\
&4\ar[d]^{\bar{b}}\ar[rd]^{\bar{c}}\\
&2\ar[d]^{b}&3\ar[ld]^{(c,-1)}\\
&4}},\qquad
G:=\vcenter{\xymatrix{
&2\ar[d]^b\\
&4\ar[ld]_{\bar{a}}\ar[rd]^{\bar{c}}\\
1\ar[rd]_a&&3\ar[ld]^{(c,-1)}\\
&4
}},\qquad
H:=\vcenter{\xymatrix{
&&3\ar[ld]^{c}\\
&4\ar[ld]_{\bar{a}}\ar[d]^{\bar{b}}\\
1\ar[rd]_a&2\ar[d]^{(b,-1)}\\
&4}},
\]
we obtain
\begin{align*}
\del_{M(0)}      &= \del_{M(\lam)}+\del_H, &\del_A =&\del_R+\del_F,\\
\del_{M(-1)}     &= \del_{M(\lam)}+\del_F, &\del_B =&\del_R+\del_G,\\
\del_{M(\infty)} &= \del_{M(\lam)}+\del_G, &\del_C =&\del_R+\del_H.\\
\end{align*}
These equalities follow from simple calculations of Euler characteristics 
of varieties of composition series. 
For example, to see the first equality
one should observe that for $M(0)$ there are two types of composition
series: Those with top $S_4$ may be identified with the composition series
of $M(\lam)$ for $\lam$ ``generic''. The remaining composition series have
top $S_3$ and may be identified with the composition series of $H$.
Our claim follows, see~\ref{ssec:repeat}.

So eventually
\[
\del_T\cdot\del_{S_4}= \del_{M(\lam)}+\del_R+\del_F+\del_G+\del_H
\]
as an expansion in the dual semicanonical basis.

%%%%%%%%%%%%%%%%%%%%%%%%%%%%%%%%%%%%%%%%%%%%%%%%%%%%%%%%%%%%%%%%%%%%%%%%%%%%%%%
\section{Recollections}

%%%%%%%%%%%%%%%%%%
\subsection{2-Calabi-Yau algebras and categories} \label{rec:2cy}
%%%%%%%%%%%%%%%%%%
Let $\ka$ be a field.
A triangulated $\ka$-category $\cT$ with suspension functor $\Sig$ is called
a {\em d-Calabi-Yau category} if all its homomorphism spaces are 
finite-dimensional over $\ka$, and there exists a functorial isomorphism
\[
\cT(x,y) \ra \Hom_\ka(\cT(y,\Sig^d x),\ka)
\]
for all $x,y\in\cT$.

It is well-known, that for a preprojective algebra $\La$  of Dynkin type 
the stable module category $\La\text{-}\smd$ with the suspension functor
$\Ome^{-1}$ is a 2-Calabi-Yau category, see for example~\cite[\S 7.5]{GLS2}.
In this case equation~\eqref{eqn:CY-Func} follows
via the natural isomorphism~\cite[\S 6]{Hel}
$\Ext^1_\La(x,y)\ra\sthom_\La(x,\Ome^{-1}y)$
from the 2-Calabi-Yau property of $\La\text{-}\smd$.

Recently it was announced~\cite{ChR} (see also~\ref{rem:extsym}~(1))
that a (connected) preprojective algebra $\La$  
which is not of Dynkin type is a {\em 2-Calabi-Yau algebra}, so by definition 
$\De_f(\La)$, the full subcategory of the bounded derived category of $\La$
of complexes with finite-dimensional nilpotent cohomology, is a
2-Calabi-Yau category. 
In this case  equation~\eqref{eqn:CY-Func} 
follows from the 2-Calabi-Yau property of $\De_f(\La)$
via the natural isomorphism $\Ext^1_\La(x,y)\ra \De_f(x,y[1])$ for
nilpotent $\La$-modules $x,y$.

%%%%%%%%%%%%%%
\subsection{Extensions} \label{pushouts}
%%%%%%%%%%%%%%
If $\eps\df 0\ra y''\ra y\ra y'\ra 0$ is a short exact sequence of
$\La$-modules, we write $[\eps]$ for its class in $\Ext^1_\La(y',y'')$.
Let us recall the functorial behavior of $\Ext^1_\La$ in terms of
short exact sequences. 
If $\rho\in\Hom_\La(x',y')$ and 
$\lam\in\Hom_\La(y'',z'')$ then $[\eps]\circ\rho\in\Ext^1_\La(x',y'')$ is
represented precisely by a short exact sequence 
$0\ra y''\ra y\ra x'\ra 0$ which is obtained from $\eps$ as the pullback
along $\rho$. 
Similarly, $\lam\circ[\eps]\in\Ext^1_\La(y',z'')$ is
represented precisely by a short exact sequence
$0\ra z''\ra z\ra y'\ra 0$ which is obtained from $\eps$ as the pushout
of $\eps$ along $\lam$. 
\[
\def\objectstyle{\scriptstyle}
\def\labelstyle{\scriptstyle}
\xymatrix @-1.2pc{
&{\textstyle\eps^\lam\df}  &0\ar[rr] &&{z''}\ar[rr]&& z \ar[rr]&& {y'}\ar[rr]&&0\\
{\textstyle\eps\df} &0\ar[rr]&&{y''}\ar[ru]^{\lam}\ar[rr]&& y\ar[ru]\ar[rr]&& 
{y'}\ar@{=}[ru]\ar[rr]&&0\\
&{\textstyle\eps_\rho^\lam\df}  &0\ar'[r][rr] &&{z''}\ar@{=}'[u][uu]\ar'[r][rr]&& e\ar'[u][uu]\ar'[r][rr] && 
{x''}\ar'[u][uu]_{\rho}\ar[rr]&&0\\
{\textstyle\eps_\rho\df}&0\ar[rr]&& {y''}\ar@{=}[uu]\ar[rr]\ar[ru]^{\lam} && x\ar[uu]\ar[rr]\ar[ru]&& 
x'\ar[uu]_(.7){\rho}\ar@{=}[ru]\ar[rr]&&0
}
\]
So we have in the above diagram $[\eps^\lam]=\lam\circ[\eps]$ and
$[\eps_\rho]=[\eps]\circ\rho$ and the associativity
$(\lam\circ[\eps])\circ\rho=[\eps^\lam_\rho]=\lam\circ([\eps]\circ\rho)$.
%%%%%%%%%%%%%
\subsection{Bilinear forms and orthogonality}\label{bilinear}
%%%%%%%%%%%%%
Let $\phi\df U \times U' \to \C$
be a bilinear form.
For subspaces $L \subseteq U$ and 
$L' \subseteq U'$ define
\begin{align*}
L^\perp &= \{ u' \in U' \mid \phi(u,u')=0 
\text{ for all }
u \in L \},\\
{^\perp}L' &=  \{ u \in U \mid 
\phi(u,u')= 0 \text{ for all }
u' \in L' \}.
\end{align*}
We call $u \in U$ and $u' \in U'$ {\it orthogonal} (with respect to
$\phi$) if $\phi(u,u') = 0$.
If $\phi$ is non-degenerate, then
${^\perp}(L^\perp) = L$ and $({^\perp}L')^\perp = L'$, and 
\[
\dim L + \dim L^\perp = \dim L' + \dim {^\perp}L' = \dim U = \dim U'.
\]

%%%%%%%%%%%%%%%%%%%%%%%%%
%\subsection{}\label{par2}
%%%%%%%%%%%%%%%%%%%%%%%%%
Let
$\phi_V\df V \times V' \to \C$ and 
$\phi_W\df W \times W' \to \C$
be non-degenerate bilinear forms, and let 
$f \in \Hom_\C(V,W)$.
Then we can identify $f$ with the dual $f'^*\df V'^* \to W'^*$ of a map
$f' \in \Hom_\C(W',V')$ if and only if
\begin{equation}\label{pairing}
\phi_V(v,f'(w')) = \phi_W(f(v),w')
\end{equation}
for all $v \in V$ and $w' \in W'$.
Here we use the isomorphisms $\widetilde{\phi}_V\df V \to V'^*$ and 
$\widetilde{\phi}_W\df W \to W'^*$ defined by
$v \mapsto \phi_V(v,-)$ and $w \mapsto \phi_W(w,-)$, respectively.
Assume that Equation (\ref{pairing}) holds.
Thus
we get a commutative diagram
\[
\xymatrix{
V \ar[r]^f \ar[d]_{\widetilde{\phi}_V} & W \ar[d]^{\widetilde{\phi}_W}\\
V'^* \ar[r]^{f'^*} & W'^*
}
\;\;\;\;\;\;\;\;\;\;
\xymatrix{
v \ar@{|->}[rr]^{f} \ar@{|->}[d]_{\widetilde{\phi}_V} && 
f(v) \ar@{|->}[d]^{\widetilde{\phi}_W} \\
\phi_V(v,-) \ar@{|->}[r]^{f'^*} & \phi_V(v,f'(-)) \ar@{=}[r] &
\phi_W(f(v),-).
}
\]
The proof of the following Lemma is an easy exercise.
\begin{Lem}\label{orthogcor2}
$\Ker(f') = \Ima(f)^\perp$.
\end{Lem}

%%%%%%%%%%%%%
\subsection{Euler characteristics}
%%%%%%%%%%%%%
For a complex algebraic variety $V$, let ${\mathcal C}(V)$ denote 
the abelian group of constructible functions on $V$ with respect 
to the Zariski topology. 

If $\pi\df V \to W$ is a morphism of complex varieties and 
$f\in{\mathcal C}(V)$, we define a function $\pi_*f$ on $W$
by
\[
(\pi_*f)(w) = \int_{\pi^{-1}(w)} f = \sum_{a \in \C} a \, 
\chi(f^{-1}(a) \cap \pi^{-1}(w)),\qquad (w\in W).
\] 
Then it is known that $\pi_*f$ is constructible.
Moreover, for morphisms 
$\pi\df V \to W$ and $\theta\df W \to U$ of complex varieties we have
\[
(\theta \circ \pi)_* = \theta_* \circ \pi_*. 
\]
(The case of compact complex algebraic varieties is discussed in
\cite[Proposition 1]{MP}, the general case can be found in 
\cite[Proposition 4.1.31]{D}.)
As a particular case, we note  the following result:
\begin{Prop}\label{eulerfibres}
Let $\pi\df V \to W$ be a morphism of complex varieties such that
there exists some $c \in \Z$ with 
$
\chi(\pi^{-1}(w)) = c
$
for all $w \in \Ima(\pi)$.
Then 
\[
\chi(V) = c\, \chi(\Ima(\pi)).
\]
In particular, if $\pi\df V \to W$ is a morphism of complex varieties such that
for all $w \in \Ima{\pi}$ the fiber
$
\pi^{-1}(w)
$
is isomorphic to an affine space, then
$\chi(V) = \chi(\Ima(\pi))$.
\end{Prop}

%%%%%%%%%%%%%%%%%%%%%%%%%%%%%%%%%%%%%%%%%%%%%%%%%%%%%%%%
\section{Ext-Symmetry for preprojective algebras} \label{sec:extsym}
The main goal of this section is to provide a direct proof of the following
result which is crucial for this paper. Otherwise it is independent of
the main body of the paper.

\begin{Thm} \label{Th:extsym}
For a quiver $Q$ without loops
let $\La$ be the associated
preprojective algebra over a field $\ka$.
Then for all finite-dimensional $\La$-modules $M$ and $N$ there is a functorial
isomorphism 
$$
\phi_{M,N}\df\Ext^1_\La(M,N)\ra D\Ext^1_\La(N,M).
$$
\end{Thm}

\subsection{Preliminaries}
Let $\ka$ be a field, and let
$Q=(Q_0,Q_1,s,e)$ be a finite quiver without loops.
Here $Q_0$ denotes the set of vertices, $Q_1$ is
the set of arrows of $Q$, and $s,e\df Q_1 \to Q_0$ are maps.
An arrow $\alp \in Q_1$ {\it starts} in a vertex $s\alp = s(\alp)$ 
and {\it ends} in a vertex
$e\alp = e(\alp)$.

By $\ol{Q}$ we denote the {\it double quiver} of $Q$, so $\ol{Q}_0=Q_0$,
$\ol{Q}_1=Q_1\cup \{\ol{\alp}\mid \alp\in Q_1\}$, and we extend the
maps $s,e$ to maps $s,e\df \ol{Q}_1 \to Q_0$ 
by
$s\ol{\alp}:=e\alp$, $e\ol{\alp}:=s\alp$ for all $\alp\in Q_1$. 
It will be convenient to consider $\ol{?}$ as an involution on
$\ol{Q}_1$ with 
\[
\ol{\bet} := \begin{cases} \ol{\bet} &\text{if } \bet\in Q_1,\\
                        \alp      &\text{if } \bet=\ol{\alp} \text{ for some }
                        \alp\in Q_1.\end{cases}
\]
Moreover, for all $\beta \in \ol{Q}_1$ we define
\[
\abs{\bet} := \begin{cases} 0 &\text{if } \bet\in Q_1,\\
                            1 &\text{else}.
              \end{cases}
\]
We consider the preprojective algebra
\[
\La =\ka\ol{Q}/\bil{\rho_i\mid_{i\in Q_0}}\quad\text{where }
\rho_i=\sum_{\substack{\bet\in\ol{Q}_1\\s\bet=i}} 
(-1)^{\abs{\bet}}\ol{\bet}\bet.
\]
Note that $\La$ is a quadratic, possibly infinite-dimensional quiver algebra, 
in any case $\La$ is augmented over  $S:=\ka^{\times Q_0}$
({\em i.e.} we have $\ka$-algebra homomorphisms $S\ra\La\ra S$ whose 
composition is the identity on $S$).
Note that $S$ is a commutative semisimple $\ka$-algebra
which has a natural basis consisting of primitive orthogonal idempotents
$\{e_i \mid i \in Q_0\}$. 
In what follows, all undecorated tensor products are meant over $S$ .

The proof  relies on the following well-known result which holds
{\em mutatis mutandis} more generally for any quadratic quiver algebra.
It follows for example from  the proof of~\cite[Theorem 3.15]{BrBuKi02}.

\begin{Lem} \label{Lem:extsym}
Let $\La$ be the preprojective algebra as defined above.
Then
\[
P^\bullet\df
\La\otimes R\otimes\La\xra{d^1}\La\otimes A\otimes\La\xra{d^0}
\La\otimes S\otimes\La
\]
is the beginning of a projective bimodule resolution of $\La$, 
where $S=\oplus_{i\in Q_0}\ka e_i$, $A=\oplus_{\bet\in\ol{Q}_1}\ka\bet$ 
and $R=\oplus_{i\in Q_0}\ka\rho_i$ are $S\mc S$-bimodules in the obvious way, 
and
\begin{align*}
d^0=&\sum_{\bet\in\ol{Q}_1}
(\bet\otimes e_{s\bet}\otimes e_{s\bet}-e_{e\bet}\otimes e_{e\bet}\otimes\bet)
\otimes_{S\mc S}\bet^*,\\
d^1 = &\sum_{\bet\in\ol{Q}_1} (-1)^{\abs{\bet}}
(\ol{\bet}\otimes\bet\otimes e_{s\bet}+e_{s\bet}\otimes\ol{\bet}\otimes\bet)
\otimes_{S\mc S}\rho_{s\bet}^*.
\end{align*}
\end{Lem}

In the statement of the Lemma we denote for example by
$(\bet^*)_{\bet\in\ol{Q}_1}$ the dual basis for the $S\mc S$-bimodule
$DA=\Hom_\ka(A,\ka)$. 
In this case we have 
$e_i\bet^* e_j=\delta_{i,s\bet}\delta_{j,e\bet}\bet^*$.

%%%%%%%%%%%%%%%%%%%%%%%%%%%%%%%%%%%%%%%%%%%%%%%%%%%%%%%%%%%%%%%%%%%%%%%%%%%%%
\subsection{Proof of Theorem~\ref{Th:extsym}} \label{ssec:extsym2}
For finite-dimensional $\La$-modules $M$ and $N$ it is
easy to determine the complex $\Hom_\La(P^\bullet\otimes_\La M,N)$:
\begin{equation}\label{eq:Cpx1}
\bigoplus_{i\in Q_0} \Hom_\ka(M(i),N(i))\xra{d^0_{M,N}}
\bigoplus_{\bet\in \ol{Q}_1}\Hom_\ka(M(s\bet), N(e\bet))\xra{d^1_{M,N}}
\bigoplus_{i\in Q_0} \Hom_\ka(M(i),N(i))
\end{equation}
where $d^0_{M,N}=\Hom_\La(d^0\otimes_\La M,N)$ is given by
\[
(f_i)_{i\in Q_0} \mapsto (N(\bet)f_{s\bet}-f_{e\bet}M(\bet))_{\bet\in\ol{Q}_1}
\]
and $d^1_{M,N}=\Hom_\La(d^1\otimes_\La M,N)$ is given by
\[
(g_\bet)_{\bet\in\ol{Q}_1} \mapsto
\left(\sum_{\bet\in\ol{Q}_1\df s\bet=i}
(-1)^{\abs{\bet}} (N(\ol{\bet})g_\bet + g_{\ol{\bet}}M(\bet)) 
\right)_{i\in Q_0}.
\]
Thus, we have functorially $\Ext^1_\La(M,N)=\Ker(d^1_{M,N})/\Ima(d^0_{M,N})$.

Similarly, we find $D\Hom_\La(P^\bullet\otimes_\La N,M)$, which is identified
via the trace pairings
\begin{multline*}
(\bigoplus_{i\in Q_0}\Hom_\ka(M(i),N(i)))\times 
(\bigoplus_{i\in Q_0}\Hom_\ka(N(i),M(i)))\ra\ka,\\
((\vph_i),(f_i))\mapsto \sum_{i\in Q_0} \Tr(\vph_i\circ f_i) 
\end{multline*}
and
\begin{multline*}
(\bigoplus_{\bet\in\ol{Q}_1}\Hom_\ka(M(s\bet),N(e\bet))\times
(\bigoplus_{\bet\in\ol{Q}_1}\Hom_\ka(N(s\bet),M(e\bet))\ra\ka,\\
((\eps_\bet),(g_\bet))\mapsto \sum_{\bet\in\ol{Q}_1} 
\Tr(\eps_{\ol{\bet}}\circ g_\bet)
\end{multline*}
to
\begin{equation}\label{eq:Cpx2}
\bigoplus_{i\in Q_0} \Hom_\ka(M(i),N(i))\xra{d^{1,*}_{N,M}}
\bigoplus_{\bet\in \ol{Q}_1}\Hom_\ka(M(s\bet), N(e\bet))\xra{d^{0,*}_{N,M}}
\bigoplus_{i\in Q_0} \Hom_\ka(M(i),N(i))
\end{equation}
with
$d^{1,*}_{N,M}=D\Hom_\La(d^1\otimes N,M)$ given by
\[
(\vph_i)_{i\in Q_0}\mapsto 
((-1)^{\abs{\bet}}(N(\bet)\vph_{s\bet}-\vph_{e\bet}M(\bet)))_{\bet\in\ol{Q}_1}
\]
and
$d^{0,*}_{N,M}=D\Hom_\La(d^0\otimes N,M)$ given by
\[
(\eps_\bet)_{\bet\in\ol{Q}_1}\mapsto
\left(
\sum_{\bet\in\ol{Q}_1\df e\bet=i}\eps_\bet M(\ol{\bet}) -
\sum_{\bet\in\ol{Q}_1\df s\bet=i} N(\ol{\bet})\eps_\bet
\right)_{i\in Q_0}.
\]
In fact, we have by definition
\begin{align*}
d^{1,*}_{N,M}((\vph_i)_{i\in Q_0})((g_\bet)_{\bet\in\ol{Q}_1}) 
&=\sum_{i\in Q_0} \Tr(\vph_i\circ(\sum_{\substack{\bet\in\ol{Q}_1\\s\bet=i}}
(-1)^{\abs{\bet}}(M(\ol{\bet})g_\bet+g_{\ol{\bet}}N(\bet))))\\
&=\sum_{\bet\in\ol{Q}_1}(-1)^{\abs{\bet}}
(\Tr(\vph_{e\ol{\bet}}M(\ol{\bet})g_\bet)
+\Tr(g_{\ol{\bet}}N(\bet)\vph_{s\bet}))\\
 &=\sum_{\bet\in\ol{Q}_1}(-1)^{\abs{\bet}}
\Tr((N(\bet)\vph_{s\bet} - \vph_{e\bet}M(\bet))\circ g_{\ol{\bet}})\\
\intertext{and}
d^{0,*}_{N,M}((\eps_\bet)_{\bet\in\ol{Q}_1})((f_i)_{i\in Q_0}) &=
\sum_{\bet\in\ol{Q}_1} \Tr(\eps_{\ol{\bet}}\circ(M(\bet)f_{s\bet}-f_{e\bet}N(\bet)))\\
&=
\sum_{i\in Q_0}
\Tr((\sum_{\substack{\bet\in\ol{Q}_1\\e\bet=i}}\eps_\bet M(\ol{\bet})-
\sum_{\substack{\bet\in\ol{Q}_1\\s\bet=i}} N(\ol{\bet})\eps_\bet)\circ\vph_i).
\\
\end{align*}
Thus, we have functorially 
$D\Ext^1_\La(N,M)=\Ker(d^{0,*}_{N,M})/\Ima(d^{1,*}_{N,M})$.
The complexes \eqref{eq:Cpx1} and \eqref{eq:Cpx2} are isomorphic:
\[
\xymatrix@-0.4pc{
\oplus_{i\in Q_0}\Hom_\ka(M(i),N(i))\ar[r]^{d^0_{M,N}}\ar[d]^{\id}&
\oplus_{\bet\in\ol{Q}_1}\Hom_\ka(M(s\bet),N(e\bet))
\ar[r]^{d^1_{M,N}}\ar[d]^{\Theta_{M,N}}& 
\oplus_{i\in Q_0}\Hom_\ka(M(i),N(i))\ar[d]^{\id}
\\
\oplus_{i\in Q_0}\Hom_\ka(M(i),N(i))\ar[r]^{d^{1,*}_{N,M}}&
\oplus_{\bet\in\ol{Q}_1}\Hom_\ka(M(s\bet),N(e\bet))\ar[r]^{d^{0,*}_{N,M}}&
\oplus_{i\in Q_0}\Hom_\ka(M(i),N(i))}
\]
with $\Theta_{M,N}((g_\bet)_{\bet\in\ol{Q}_1})=
((-1)^{\abs{\bet}}g_\bet)_{\bet\in\ol{Q}_1}$.
This finishes the proof of Theorem~\ref{Th:extsym}.

\subsection{Remarks} \label{rem:extsym}
(1) Assume $Q$ is connected.
Then it is known that in the situation of Lemma~\ref{Lem:extsym}
we have an isomorphism of bimodules:
\[
\Ker{d^1} \cong\begin{cases}
D\La &\text{ if } Q \text{ is a Dynkin quiver,}\\
\;0     &\text{ else.}
\end{cases}
\]
This  follows for example from~\cite[Thm.~4.8. \& Thm.~4.9]{BrBuKi02}
in the first case and from~\cite[Prop.~4.2]{BrBuKi02} together 
with~\cite[Thm.~9.2]{BuKi99} in the second case.

Thus, $P^\bullet$ is a projective bimodule resolution of $\La$, if
$Q$ is not a Dynkin quiver.
A similar calculation as in~\ref{ssec:extsym2} shows that
$P^\bullet$ is self-dual in the sense of Bocklandt~\cite[\S 4.1]{Boc},
thus the bounded derived category $\mathcal{D}^b(\La$-$\md)$ is
2-Calabi-Yau~\cite[Thm.~4.2]{Boc}. 

(2) We leave it as an (easy) exercise to derive from the calculations 
in Section~\ref{ssec:extsym2} 
Crawley-Boevey's formula~\cite[Lemma~1]{cb00},
and in the non-Dynkin case the equation
\[
\dim\Hom_\La(M,N)-\dim\Ext^1_\La(M,N)+\dim\Ext^2_\La(M,N)=
(\dimv M,\dimv N)_Q,
\]
where $(-,-)_Q$ is the symmetric quadratic form associated to $Q$.

%%%%%%%%%%%%%%%%%%%%%%%%%%%%%%%%%%%%%%%%%%%%%%%%%%%%%%%%%%%%%%%%%%%%%%

\vspace{1.5cm}
%%%%%%%%%%%%%%%%%%%%%%%%%%%%%%%%%%%%%%%
{\parindent0cm \bf Acknowledgements.}\,
%%%%%%%%%%%%%%%%%%%%%%%%%%%%%%%%%%%%%%%
We like to thank Bernhard Keller and William Crawley-Boevey 
for interesting and helpful discussions.

%%%%%%%%%%%%%%%%%%%%%%%%%%%%%%%%%%%%%%%%%%%%%%%%%%%%%%%%%


\begin{thebibliography}{CH}

%%%%%%%%%%%%%%%%%%%%%%%%%%%%%%%%%%%%%%%%%%%%%%%%%%%%%%%%%

\bibitem{Boc}
{\it R. Bocklandt},
Graded Calabi Yau Algebras of dimension 3.
arXiv:math.RA/0603558.

\bibitem{Bonga}
{\it K. Bongartz},
Some geometric aspects of representation theory.
Algebras and modules, I (Trondheim, 1996),  1--27, CMS Conf. Proc., 23, 
Amer. Math. Soc., Providence, RI, 1998.

\bibitem{BrBuKi02}
{\it S.~Brenner, M.C.R. Butler,  A.D. King}, 
Periodic algebras which are almost {K}oszul.
Algebr. Represent. Theory 5 (2002), no.~4, 331--367.

\bibitem{BuKi99}
{\it M.~C.~R. Butler  A.~D. King}, 
Minimal resolutions of algebras.
J. Algebra 212 (1999), no.~1, 323--362.

\bibitem{ChR}
{\it J. Chuang, R. Rouquier},
Paper in preparation.

\bibitem{CK}
{\it P. Caldero, B. Keller},
From triangulated categories to cluster algebras.
arXiv:math.RT/0506018.

\bibitem{cb00}
{\it W.~Crawley-Boevey}, 
On the exceptional fibres of {K}leinian singularities. 
Amer. J. Math. 122 (2000), no.~5, 1027--1037.

\bibitem{D}
{\it A. Dimca},
Sheaves in topology.
Universitext. Springer-Verlag, Berlin, 2004. xvi+236 pp. 

\bibitem{GLS}
{\it Ch. Gei{\ss}, B. Leclerc, J. Schr\"oer},
Semicanonical bases and preprojective algebras.
Ann. Sci. {\'E}cole Norm. Sup. (4) 38 (2005),  no. 2, 193--253.

\bibitem{GLS2}
{\it Ch. Gei{\ss}, B. Leclerc, J. Schr\"oer},
Auslander algebras and initial seeds for cluster algebras.
arXiv:math.RT/0506405, 
to appear in J. London Math.~Soc.

\bibitem{GLS3}
{\it Ch. Gei{\ss}, B. Leclerc, J. Schr\"oer},
Rigid modules over preprojective algebras.
Invent. Math. 165 (2006), no. 3, 589--632.

\bibitem{Hel}
{\it A. Heller},
The loop-space functor in homological algebra.  
Trans. Amer. Math. Soc.  96 (1960), 382--394.

\bibitem{Lu}
{\it G. Lusztig},
Quivers, perverse sheaves, and quantized enveloping algebras. 
J. Amer. Math. Soc. 4 (1991), no. 2, 365--421.

\bibitem{Lu2}
{\it G. Lusztig},
Semicanonical bases arising from enveloping algebras.  
Adv. Math.  151  (2000),  no. 2, 129--139.

\bibitem{MP}
{\it R.D. MacPherson},
Chern classes for singular algebraic varieties.
Ann. of Math.(2) 100 (1974), 423--432.

\bibitem{Ri1}
{\it C. M. Ringel},
The preprojective algebra of a quiver. 
In: Algebras and modules II (Geiranger, 1966),
467--480, CMS Conf. Proc. 24, AMS 1998.

\bibitem{Serre}
{\it J.P. Serre},
Espaces fibr{\'e}s alg{\'e}briques.
Seminaire C. Chevalley 1958, 1--36.


\end{thebibliography}
\end{document}